\documentclass[aos,MSNbibl,seceqn,nameyear,dvips]{arximspdf}

\doi{10.1214/12-AOS1082} 
\volume{41}
\issue{1}
\pubyear{2013}
\firstpage{296}
\lastpage{322}

\makeatletter

\newtheorem{theorem}{Theorem}[section]
\newtheorem{proposition}{Proposition}[section]
\newtheorem{lemma}{Lemma}[section]

\newproclaim{remark}{Remark}[section]

\def\b{\mathbf{b}}
\def\Ga{\bolds{\Gamma}}
\def\R{\mathbf{R}}
\def\D{\mathbf{D}}
\def\A{\mathbf{A}}
\def\W{\mathbf{W}}
\def\X{\mathbf{X}}
\def\I{\mathbf{I}}
\def\Y{\mathbf{Y}}
\def\Z{\mathbf{Z}}
\def\S{\bolds{\Sigma}}
\def\m{\bolds{\mu}}

\def\pr{\mathsf{P}}
\def\ep{\mathsf{E}}
\def\Cov{\operatorname{\mathsf{Cov}}}
\def\Var{\operatorname{\mathsf{Var}}}

\makeatother

\begin{document}
\begin{frontmatter}

\title{A Cram\'{e}r moderate deviation theorem for Hotelling's $T^{2}$-statistic with applications to~global~tests}
\runtitle{Moderate deviation for Hotelling's $T^{2}$-statistic}

\begin{aug}
\author[A]{\fnms{Weidong} \snm{Liu}\corref{}\thanksref{t2}\ead[label=e1]{liuweidong99@gmail.com}}
\and
\author[B]{\fnms{Qi-Man} \snm{Shao}\thanksref{t3}\ead[label=e3]{qmshao@cuhk.edu.hk}}
\runauthor{W. Liu and Q.-M. Shao}
\affiliation{Shanghai Jiao Tong University and Chinese University of
Hong Kong}
\address[A]{Department of Mathematics \\
Institute of Natural Sciences \\
Shanghai Jiao Tong University \\
Shanghai, P.R. China\\
\printead{e1}}
\address[B]{Department of Statistics\\
Chinese University of Hong Kong\\
Shatin, N.T.\\
Hong Kong, P.R. China\\
\printead{e3}} 
\end{aug}

\thankstext{t2}{Supported by NSFC, Grant 11201298, the Program for
Professor of
Special Appointment (Eastern Scholar) at Shanghai Institutions of
Higher Learning, the Foundation for the Author of National Excellent
Doctoral Dissertation of PR China and the startup fund from SJTU.}

\thankstext{t3}{Supported in part by Hong Kong RGC UST603710 and CUHK2130344.}

\received{\smonth{8} \syear{2012}}
\revised{\smonth{12} \syear{2012}}

%
\begin{abstract}
A Cram\'{e}r moderate deviation theorem for Hotelling's
$T^{2}$-statistic is proved under a finite $(3+\delta)$th moment. The
result is applied to large scale tests on the equality of
mean\vspace*{1pt} vectors and is shown that the number of tests can be
as large as $e^{o(n^{1/3})}$ before the chi-squared distribution
calibration becomes inaccurate. As an application of the moderate
deviation results, a global test on the equality of $m$ mean vectors
based on the maximum of Hotelling's $T^{2}$-statistics is developed and
its asymptotic null distribution is shown to be an extreme value type I
distribution. A novel intermediate approximation to the null
distribution is proposed to improve the slow convergence rate of the
extreme distribution approximation. Numerical studies show that the new
test procedure works well even for a small sample size and performs
favorably in analyzing a breast cancer dataset.
\end{abstract}

%

\pdfkeywords{62E20, 62H15, 60F10, Cramer moderate deviation,
Hotelling's T2-statistic,
global tests,
simultaneous hypothesis tests, FDR,
brain structure, gene selection}

\end{frontmatter}

\section{Introduction}\label{sec1}

Consider the following $m$ simultaneous tests:
%
\begin{equation}
\label{a0} H_{0i}\dvtx  \m_{1i}=\m_{2i} \quad\mbox{versus}\quad
H_{1i}\dvtx  \m_{1i}\neq\m_{2i}
\end{equation}
for $1\leq i\leq m$, where $\m_{1i}$ and $\m_{2i}$ are
$d_{i}\geq1$-dimensional mean vectors, and $d_{i}$ are uniformly
bounded. When $d_{i}=1$, the multiple testing problem (\ref{a0}) has
been extensively studied. A common statistical method is the two sample
$t$-test together with multiple comparison procedure by controlling the
familywise error rate (FWER) or the false\vadjust{\goodbreak} discovery rate (FDR). The
theoretical justification of this method can be found in
\citet{FanHalYao07}. Although not much attention has been paid to
the multivariate case $d_{i}>1$, (\ref{a0}) has arisen from several
important applications including shape analysis of brain structures and
gene selection.

\begin{itemize}
%
\item
\textit{Shape analysis of brain structures.} There is a growing interest
in statistical shape analysis within the neuroimaging
community; see \citet{StyOguXu}, \citet{Zhaetal}, \citet{Geretal09}.
\citet{StyOguXu} developed a widely-used software
to locate significant shape changes between healthy and pathological
brain structures. The final and most important step in \citet{StyOguXu}
procedure is the simultaneous testing of (\ref{a0}) with $\m_{1i}$
and $\m_{2i}$ being mean vectors of
3 coordinates of surface points. The number of tests $m$ can be {
hundreds or even thousands} and $d_{i}=3$ for all $i$.
In \citet{StyOguXu},
two sample Hotelling's $T^2$-statistics $T^{2}_{ni}$ were used for each
$H_{0i}$ and Benjamini--Hochberg procedure
was used to control the FDR.

\item
\textit{Gene selection.} In the breast cancer dataset analyzed by
\citet{Maretal05}, every gene corresponds to a two to
six-dimensional vector that represents the DNA methylation status of
CpG sites. Dimension $d_{i}$ is between 2 to 6. In
\citet{Maretal05}, two sample Hotelling's $T^2$-statistics and
Benjamini--Hochberg FDR correction were used to identify the
significantly different genes between two patient groups.
\end{itemize}

It is well known that Hotelling's $T^2$-statistic is asymptotically
chi-squared distributed when the
underlying distribution has a finite second moment. This provides a
natural way to
estimate $p$-values.
In the ``large $m$ small $n$'' statistical analysis, the true $p$-values
are typically small, of order $O(1/m)$
in FDR procedure. A~basic question is:

\begin{quote}
{with how many tests can the
chi-squared distribution calibration be applied before the tests become
inaccurate}?
\end{quote}

As discussed in \citet{FanHalYao07} and \citet{LiuSha10}, the
question can be
answered with Cram\'{e}r-type moderate deviation results.
The moderate deviation behavior 
for $t$-statistic is now well-understood, however, a Cram\'{e}r type
moderate deviation theorem for Hotelling's $T^2$-statistic is still not
available.
The main purpose of this paper is to establish the moderate deviation
theorem for
Hotelling's $T^2$-statistic (one-sample and two-sample). We shall prove
that under a finite $(3+\delta)$th moment,
Hotelling's $T^2$-statistic $T^{2}_{n}$ satisfies
\[
\frac{\pr(T^{2}_{n}\geq x^{2} )}{\pr(\chi^{2}(d)\geq
x^{2} )}\rightarrow1
\]
uniformly for $x\in[0, o(n^{1/6}))$.
Consequently, the number of tests can be as large as $e^{o(n^{1/3})}$
before the chi-squared distribution
calibration becomes inaccurate; see (\ref{th1b}).

As an application of the moderate deviation result, we consider the
global testing
%
\begin{eqnarray}
\label{glob} &&H_{0}\dvtx  \m_{1i}=\m_{2i} \qquad\mbox{for
all $1\leq i\leq m$}\quad \mbox{against}
\nonumber\\[-8pt]\\[-8pt]
&& H_{1}\dvtx  \m_{1i}\neq\m_{2i} \qquad\mbox{for some
$i$.}
\nonumber
\end{eqnarray}
In shape analysis of brain structures with $d_{i}=3$, the global test
(\ref{glob}) is often used to
determinate whether two brain shapes between two groups of subjects are
different or not; see \citet{CaoWor99},
\citet{TayWor08}. In gene selection [\citet{Maretal05}],
(\ref{glob}) has been used to test whether the endocrine therapy is
effective on DNA methylation status.
Here we are particularly interested in the alternative hypothesis that
the locations where $\m_{1i}\neq\m_{2i}$ are sparse.
For example, in the brain structures, the shape differences are
commonly assumed to be confined to
a small number of isolated regions inside the whole brain. In this
paper, we
shall propose a testing procedure based on the maximum of Hotelling's
$T^{2}$-statistics.
The proposed test procedure shares several advantages. {It is quite
robust} to the tails of
the underlying distribution and the dependence structure.
It converges to the given {significance} level with a rate of $\sqrt
{(\log m)^{5}/n}$.
A numerical study shows that the test procedure works quite well even
for small samples. 

The rest of our paper is organized {as follows}. In Section~\ref{sec2}, we state
Cram\'{e}r moderate deviation results for
Hotelling's $T^2$-statistic.
In Section~\ref{sec3}, we introduce our test procedure for the global test (\ref{glob}).
Theoretical results of the robustness on the tails and dependence
structures are given. The power of the test procedure is also
investigated.
A numerical study is carried out in Section~\ref{sec4}, in which we compare our
test procedure to some existing test procedures.
The proofs of the main results are postponed to Section~\ref{sec5}.



\section{A Cram\'{e}r type moderate deviation theorem for Hotelling's
$T^{2}$-statistic}\label{sec2}

The properties of Hotelling's $T^{2}$-statistic under normality are well known [\citet{And03}]. Large and
moderate deviations (logarithm of the tail probabilities) were obtained
in \citet{DS06}. In this section, we shall establish
a Cram\'er moderate deviation theorem for Hotelling's $T^2$-statistic.
For Student \mbox{$t$-statistic},
the Cram\'er moderate deviation result was first obtained by \citet{shao2}
under a finite third moment and the result was extended to
self-normalized sums of
independent random variables in \citet{jing}. We refer to
\citet{dela}
for a systematic presentation on the self-normalized limit theory and
its statistical applications.

Let $\{\X_{1},\ldots,\X_{n_{1}}\}$ and $\{\Y_{1},\ldots,\Y_{n_{2}}\}
$ be {two groups of} i.i.d. $d$-dimensional random vectors
with mean vectors $\m_{1}$ and $\m_{2}$ and covariance matrices $\S
_{1}$ and $\S_{2}$, respectively. Assume that
$\{\X_{1},\ldots,\X_{n_{1}}\}$ and $\{\Y_{1},\ldots,\Y_{n_{2}}\}$
are independent {and $\S_{1}$ and $\S_{2}$ are positive definite.} Let
\[
\bar{\X}=\frac{1}{n_{1}}\sum_{k=1}^{n_{1}}
\X_{k},\qquad \bar{\Y}=\frac{1}{n_{2}}\sum_{k=1}^{n_{2}}
\Y_{k}
\]
be the sample means and
\[
\mathbf{V}_{n1}=\frac{1}{n_{1}}\sum_{k=1}^{n_{1}}(
\X_{k}-\bar{\X}) (\X_{k}-\bar{\X})^{\prime},\qquad
\mathbf{V}_{n2}=\frac{1}{n_{2}}\sum_{k=1}^{n_{2}}(
\Y_{k}-\bar{\Y}) (\Y_{k}-\bar{\Y})^{\prime}
\]
be the sample covariance matrices, where for a vector $\mathbf
{a}$, $\mathbf{a}^{\prime}$ denotes its transpose.
The two sample Hotelling's $T^{2}$-statistic is then defined by
\[
T^{2}_{n}=(\bar{\X}-\bar{\Y})^{\prime} \biggl(
\frac{1}{n_{1}}\mathbf{V}_{n1}+\frac{1}{n_{2}}
\mathbf{V}_{n2} \biggr)^{-1}(\bar{\X}-\bar{\Y}).
\]
Let $n_{1}\asymp n_{2}$ denote the inequality $c_{1}\leq
n_{1}/n_{2}\leq c_{2}$ for some positive constants $c_{1}$ and $c_{2}$.
The following result gives a Cram\'{e}r type moderate deviation for
Hotelling's $T^{2}$-statistic.

\begin{theorem}\label{th1} Suppose that $n_{1}\asymp n_{2}$, $\ep\|\X
_{1}\|^{3+\delta}<\infty$ and $\ep\|\Y_{1}\|^{3+\delta}<\infty$ for
some $\delta>0$. Then, under $\m_{1}=\m_{2}$
%
\begin{equation}\label{th1a}
\frac{\pr(T^{2}_{n}\geq x^{2} )}{\pr(\chi^{2}(d)\geq
x^{2} )}\rightarrow1 \qquad\mbox{as $n\rightarrow\infty$}
\end{equation}
uniformly for $x\in[0, o(n^{1/6}))$, where $n=n_{1}+n_{2}$.
\end{theorem}

Theorem~\ref{th1} shows\vspace*{2pt} that the true distribution of $T^{2}_{n}$ can
be well approximated by $\chi^{2}(d)$ distribution
uniformly\vspace*{1pt} in the interval $[0, o(n^{1/3}))$ under the finite $(3+\delta
)$th moment. Let
$F_{n}(x)=\pr(T^{2}_{n}\geq x|\m_{1}=\m_{2} )$ and
$F(x)=\pr(\chi^{2}(d)\geq x )$.
Then,\vspace*{1pt} the true $p$-value { is} $p_{n}=F_{n}(T^{2}_{n})$ and the
estimated $p$-value is $\hat{p}_{n}=F(T^{2}_{n})$. Thus by
(\ref{th1a}),
%
\begin{equation}\label{th1b}
\biggl|\frac{\hat{p}_{n}}{p_{n}}-1 \biggr|I\bigl\{p_{n}\geq e^{-o(n^{1/3})}\bigr\}=o(1).
\end{equation}
This provides a theoretical justification of the accuracy of the
estimated $p$-values by the chi-squared distribution used in B-H FDR
correction method. We refer to \citet{FanHalYao07} and \citet{LiuSha10} for more detailed discussion on the relations between the
Cram\'{e}r type moderate deviation and the accuracy of the estimated
$p$-values used in large scale tests.

For one-sample Hotelling's $T^2$-statistic, we have a similar result.

\begin{theorem} \label{th1-2}
Suppose that $\ep\|\X_{1}\|^{3+\delta}<\infty$ for
some $\delta>0$. Then
%
\begin{equation}\label{th1-2a}
\frac{\pr( n_1 ( \bar{\X} - \m_1)^{\prime} \mathbf{V}_{n1}^{-1}
(\bar{\X}- \m_1) \geq x^{2} )} {
\pr(\chi^{2}(d)\geq x^{2} )}\rightarrow1 \qquad\mbox{as $n_{1}\rightarrow
\infty$}
\end{equation}
uniformly for $x\in[0, o(n^{1/6}_{1}))$.\vadjust{\goodbreak}
\end{theorem}

The proof of Theorem~\ref{th1-2} is completely similar to that of
Theorem~\ref{th1} and so will be omitted.

\begin{remark} As proved by \citet{shao2} and \citet{jing},
(\ref{th1a}) and (\ref{th1-2a}) hold
under finite third moments when $d=1$ and the range $[0, o(n^{1/6}))$
is the widest possible. We conjecture that (\ref{th1a}) and (\ref{th1-2a})
remain valid for $d \geq2$ under a finite third moment
and that the range $[0, o(n^{1/6}))$ is optimal.
\end{remark}

\section{Global testing}\label{sec3}

In this section, we are interested in the global testing
(\ref{glob}), that is,
\begin{eqnarray*}
&&H_{0}\dvtx  \m_{1i}=\m_{2i}\qquad\mbox{for
all $1\leq i\leq m$} \quad\mbox{against}
\\
&& H_{1}\dvtx  \m_{1i}\neq\m_{2i} \qquad\mbox{for some
$i$.}
\end{eqnarray*}
%
%
where $\m_{1i}$ and $\m_{2i}$ are $d_{i}$-dimensional mean vectors of
random vectors $\X^{i}$ and $\Y^{i}$, respectively.

Write $\mathbf{a}=(\m^{\prime}_{11},\ldots,\m^{\prime}_{1m})$ and $\b
=(\m^{\prime}_{21},\ldots,\m^{\prime}_{2m})$. Most of existing works on
the global tests are focused on the alternative that $\mathbf{a}-\b$ is
either sparse or dense. When the alternative is sparse, the commonly
used test statistic is the maximum of univariate $t$-statistics and the
higher criticism (HC$^{*}$) test procedure [\citet{DonJin04,HalJin10}]. On
the other hand, if the signals are dense, then the squared sum type
test statistics have been used [\citet{CheQin10}]. In this
section, we focus on the sparse alternative hypothesis. The main
difference between the current paper and the previous works is that the
sparse signals appear in groups and that the underlying distributions
are not necessarily normal and the components may not have an ordered
structure. For the sparse case, it has been proved in
\citet{DonJin04} that the higher criticism statistic enjoys some
optimal properties with respect to the detection region. On the other
hand, the independence between variables plays an important role in the
control of type I errors of the higher criticism statistic. The
simulation in Section~\ref{sec4} shows that HC$^{*}$ statistic may not
be robust against the dependence and may fail to control the type I
error. In contrast, our test procedure introduced below is robust to
dependence, as shown by Theorems~\ref{th3}--\ref{th5-5} and the
simulation.

Suppose that we have two groups of i.i.d. observations
\[
\mathcal{X}=\bigl\{\X^{1}_{k},\ldots,
\X^{m}_{k}; 1\leq k\leq n_{1}\bigr\}
\quad\mbox{and}\quad \mathcal{Y}=\bigl\{\Y^{1}_{k},\ldots,
\Y^{m}_{k}; 1\leq k\leq n_{2}\bigr\}
\]
with mean vectors $\{\m_{11},\ldots,\m_{1m}\}$ and
$\{\m_{21},\ldots,\m_{2m}\}$, respectively. The two groups of
observations $\mathcal{X}$ and $\mathcal{Y}$ are independent. Let
$T^{2}_{ni}$ be the\vspace*{1pt} two sample Hotelling's
$T^{2}$-statistics based on $\{\X^{i}_{k}; 1\leq k\leq n_{1}\}$ and
$\{\Y^{i}_{k}; 1\leq k\leq n_{2}\}$. We introduce our test procedure as
follows.

\textit{Case} 1. \textit{$d_{i}\equiv d$.} Let $\W_{1,k}$, $1\leq k\leq n_{1}$,
and $\W_{2,k}$, $1\leq k\leq n_{2}$ be i.i.d. multivariate normal
vectors with mean zero and covariance matrix $\I_{d}$. Let
%
\begin{equation}
\label{a10} F_{n_{1},n_{2}}(y)=\pr\bigl(T^{*2}_{n}\geq y
\bigr),
\end{equation}
where $T^{*2}_{n}$ is the two sample Hotelling's $T^2$-test statistic
based on $\{\W_{1,k}\}$ and $\{\W_{2,k}\}$.
For given $0 < \alpha< 1$, let $y_{n}(\alpha)$ satisfy
%
\begin{equation}\label{y-alpha}
\exp\bigl(-mF_{n_{1},n_{2}}\bigl(y_{n}(\alpha)\bigr) \bigr)=1-
\alpha.
\end{equation}
Note that $1-F_{n_{1},n_{2}}(y)$ is closely related to $F$
distribution. In general, we can use simulation to obtain $y_{n}(\alpha)$.
Our test procedure for (\ref{glob}) is $\Phi^{*}_{\alpha}$, where
%
\begin{equation}\label{phi-1}
\Phi^{*}_{\alpha}=I\Bigl\{\max_{1\leq i\leq m}T^{2}_{ni}
\geq y_{n}(\alpha)\Bigr\}.
\end{equation}
The hypothesis $H_{0}$ is rejected whenever $\Phi^{*}_{\alpha}=1$.

\textit{Case} 2. \textit{$d_{i}$ may be different}. Let $F_{n_{1},n_{2},d_{i}}(y)$
be defined as in (\ref{a10}) with $d$ being replaced with $d_{i}$. Let
$G_{n_{1},n_{2},d_{i}}(y)=1-F_{n_{1},n_{2},d_{i}}(y)$. We now define
\[
\Phi^{\dag}_{\alpha}=I\Bigl\{\max_{1\leq i\leq
m}G_{n_{1},n_{2},d_{i}}
\bigl(T^{2}_{ni}\bigr)\geq g_{m}(\alpha)\Bigr\}
\]
with $g_{m}(\alpha)=1+m^{-1}\log(1-\alpha)$. The hypothesis $H_{0}$
is rejected whenever \mbox{$\Phi^{\dag}_{\alpha}=1$}.
Note that $\Phi^{\dag}_{\alpha}=\Phi^{*}_{\alpha}$ if $d_{i}\equiv d$.

\begin{remark}
By Theorem\vspace*{1pt}~\ref{th3}, $\max_{1\leq i\leq m}T^{2}_{ni}$
converges to the extreme I type distribution. It seems natural to define
the following test $\Phi_{\alpha}$:
%
\begin{equation}\label{phi-2}
\Phi_{\alpha}=I\Bigl\{\max_{1\leq i\leq m}T^{2}_{ni}
\geq2\log m+(d-2)\log\log m+q_{\alpha}\Bigr\},
\end{equation}
where 
$q_{\alpha}=
-2\log(\Gamma(d/2))-2\log\log(1-\alpha)^{-1}
$. 
The hypothesis $H_{0}$ is rejected whenever $\Phi_{\alpha}=1$.
However, it is well known that
the rate of convergence to the extreme distribution is very slow [see
\citet{LiuLinSha08}].
On the other hand, the intermediate approximation given in Theorem \ref
{th5} can substantially improve the convergence rate.
This leads to our test procedure $\Phi^{*}_{\alpha}$. Numerical
results in Section~\ref{sec4} show that $\Phi^{*}_{\alpha}$ outperforms $\Phi
_{\alpha}$ significantly
and it works well even when the sample size is small.
\end{remark}


\subsection{\texorpdfstring{The limiting distribution of $\max_{1\leq i\leq m}T^{2}_{ni}$}
{The limiting distribution of max 1<=i<=m T 2 ni}}\label{sec3.1}

In this subsection,\vspace*{1pt} we show that the type I error of $\Phi^{*}_{\alpha
}$ will converges to $\alpha$ under
some mild moment conditions and dependence structure. To this end, we
need to establish the limiting distribution of
$\max_{1\leq i\leq m}T^{2}_{ni}$ under $H_0$. Let $\S_{i}=\S_{i1}+\frac
{n_{1}}{n_{2}}\S_{i2}$,
where $\S_{i1}$ and $\S_{i2}$ are the covariance matrices of $\X^{i}$
and $\Y^{i}$, respectively.
Define
\[
\Ga_{ij}=\S_{i}^{-1/2}\biggl(\Cov\bigl(
\X^{i},\X^{j}\bigr)+\frac
{n_{1}}{n_{2}}\Cov\bigl(
\Y^{i},\Y^{j}\bigr)\biggr)\S_{j}^{-1/2}.
\]
The matrix $\Ga_{ij}$ characterizes the dependence structure between
$\{\X^{i},\Y^{i}\}$ and $\{\X^{j},\Y^{j}\}$. For example, when $n_{1}=n_{2}$
and $\S_{i1}=\S_{i2}$,
\[
\Ga_{ij} = \tfrac{1}{2}\Cov\bigl(\S^{-1/2}_{i1}
\X^{i},\S^{-1/2}_{j1}\X^{j}\bigr)+
\tfrac{1}{2}\Cov\bigl(\S^{-1/2}_{i2}\Y^{i},
\S^{-1/2}_{j2}\Y^{j}\bigr)
\]
is the sum of two matrices.
When $d=1$ and $\S_{i1}=\S_{i2}$, then $\Ga_{ij}=\rho_{ij1}$, where
$\rho_{ij1}$ is the correlation coefficient between
$\X^{i}$ and $\X^{j}$.
For $0<r<1$, let
\[
\Lambda(r)=\bigl\{1\leq i\leq m\dvtx  \|\Ga_{ij}\| \geq r \mbox{ for some $j
\not= i$}\bigr\},
\]
where $\|\cdot\|$ is the spectral norm. {$\Lambda(r)$ is a subset of
$\{1,2,\ldots,m\}$ in which $\{\X^{i},\Y^{i}\}$ can be highly
correlated with other random vectors.}
Let $\R_{1}=(r_{ij1})$ and $\R_{2}=(r_{ij2})$ be the correlation
matrices of the random vectors $((\X^{1})^{\prime},\ldots,(\X^{m})^{\prime})$
and $((\Y^{1})^{\prime},\ldots,(\Y^{m})^{\prime})$, respectively.
For some $\gamma>0$, let
\[
s_{j}(m)=\operatorname{Card}\bigl\{1\leq i\leq m\dvtx  |r_{ij1}|
\geq(\log m)^{-1-\gamma
}\mbox{ or }|r_{ij2}|\geq(\log
m)^{-1-\gamma}\bigr\}.
\]
We need the following condition on the {dependence} structure.

\begin{longlist}[(C1)]
\item[(C1)] Suppose that $\operatorname{Card}(\Lambda(r))=o(m)$ for some
$0<r<1$ and
\[
\max_{1\leq j\leq p}s_{j}(m)=O(m^{\rho})
\]
for all\vspace*{1pt} $\rho>0$. Assume that
$\min_{1\leq i\leq p}\{\lambda_{\min}(\S_{i})\}\geq\tau$
for some $\tau>0$, where $\lambda_{\min}(\S_{i})$ is the smallest
eigenvalue of $\S_{i}$.
\end{longlist}

The dependence\vspace*{1pt} condition (C1) is mild. In (C1), $o(m)$ vectors $\{\X
^{i},\Y^{i}\}$, $i\in\Lambda(r)$, can be highly correlated with
other random vectors. Every $\{\X^{i},\Y^{i}\}$ can be highly
correlated with $s_{i}(m)$ vectors and weakly correlated with the
remaining vectors. The dependence in (C1) is more general than ``clumpy
dependence'' [\citet{StoTib}] and
may be a more realistic form of dependence in DNA microarrays. See also
\citet{HalWan10} who noted that short-range dependence, and more specially,
{$k$-dependence} structure,
are often observed in DNA microarrays.


The next condition is on the moment of the underlying distributions and
the relation between the sample sizes and dimension $m$.
We assume that $m$ is a function of $n=n_{1}+n_{2}$ and $m\rightarrow
\infty$ as $n\rightarrow\infty$.

\begin{longlist}[(C2)]
\item[(C2)] Suppose that $\max_{1\leq i\leq m}\ep(\|\X^{i}\|^{3+\delta
}+\|\Y^{i}\|^{3+\delta})\leq\kappa$ for some $\kappa>0$
and $\delta>0$,
$n_{1}\asymp n_{2}$ and $\log m=o(n^{1/3})$.
\end{longlist}

\begin{theorem}\label{th3} Under $H_{0}$, $d_{i}\equiv d$, \textup{(C1)} and
\textup{(C2)}, we have as $n\rightarrow\infty$,
%
\begin{eqnarray}\label{th3a}
&&\pr\Bigl(\max_{1\leq i\leq m}T^{2}_{ni}-2\log
m+(2-d)\log\log m\leq y \Bigr)
\nonumber\\[-8pt]\\[-8pt]
&&\qquad \to \exp\biggl(-\frac{1}{\Gamma(d/2)}e^{-y/2} \biggr)\nonumber
\end{eqnarray}
for any $y\in R$.
\end{theorem}

It follows from Theorem~\ref{th1} that
\[
y_{n}(\alpha)=2\log m+(d-2)\log\log m+q_{\alpha}+o(1),
\]
which together with Theorem~\ref{th3}, yields the following theorem.

\begin{theorem}\label{th3-3} Under $H_{0}$, $d_{i}\equiv d$, \textup{(C1)} and
\textup{(C2)}, we have as $n\rightarrow\infty$,
%
\begin{equation}\label{th3-3a}
\pr\bigl(\Phi^{*}_{\alpha}=1 \bigr)\rightarrow\alpha.
\end{equation}
\end{theorem}

\begin{remark}
When $d_{i}$ are different, we have a similar result as Theorem~\ref{th3-3}.
Under $H_{0}$, (C1) and (C2), we have as $n\rightarrow\infty$,
%
\begin{equation}
\label{a15} \pr\bigl(\Phi^{\dag}_{\alpha}=1 \bigr)\rightarrow
\alpha
\end{equation}
for any $0<\alpha<1$. The proof of (\ref{a15}) is similar to that of
Theorem~\ref{th3} and hence will be omitted.
\end{remark}

As mentioned earlier, the convergence rate of (\ref{th3a}) is very slow.
In testing diagonal covariance matrix problem,
\citet{LiuLinSha08} proposed to use an intermediate approximation
and proved that the rate of convergence
can be of order of $\sqrt{(\log m)^{5}/n}$. Here we give a similar
intermediate approximation to the
distribution of $\max_{1\leq i\leq m}T^{2}_{ni}$.

Let $\Theta_{j}$ be the set of indices such that $T^{2}_{nj}$ is
independent with $(T^{2}_{ni}; i\in\Theta_{j})$ and put
$s_{j}(m)=m-\operatorname{Card}(\Theta_{j})$.

\begin{longlist}[(C3$^{*}$)]
\item[(C1$^{*}$)] Suppose that $\operatorname{Card}(\Lambda(r))=O(m^{\xi})$
for some $0<r<1$ and $0\leq\xi<1$.
Assume that $\max_{1\leq j\leq m}s_{j}(m)=O(m^{\rho})$ for some
$0<\rho< (1-r)/(1+r)$.

\item[(C2$^{*}$)] Suppose that $\max_{1\leq i\leq m}\ep(\|\X^{i}\|
^{3+\delta}+\|\Y^{i}\|^{3+\delta})\leq\kappa$ for some
$\kappa>0$ and $\delta>0$,
$c_{1}\leq n_{1}/n_{2}\leq c_{2}$ for some $c_{1}>0$ and $c_{2}>0$
and $\log m=o(n^{1/3})$.

\item[(C3$^{*}$)] Suppose that $\S_{1i}=\S_{2i}$ for $1\leq
i\leq m$. We assume that $\X^{i}$ and $\Y^{i}$ can be written as
the transforms of independent components:
\[
\X^{i}=\S^{1/2}_{1i}\Z_{1i}+
\m_{1i}\quad\mbox{and}\quad \Y^{i}=\S^{1/2}_{2i}
\Z_{2i}+\m_{2i},
\]
where $\ep\Z_{1i}=0$, $\Cov(\Z_{1i})=\I$ and $\ep\Z_{2i}=0$,
$\Cov(\Z_{2i})=\I$ and the components in $\Z_{1i}$ and $\Z_{2i}$
are independent.
\end{longlist}

(C1$^{*}$) is a technical condition. It allows $T^{2}_{nj}$ be
dependent with $O(m^{\rho})$ others. By (C1$^{*}$), we can use the
Poisson approximation in \citet{ArrGolGor89}.
(C3$^{*}$) is also required for technical reason. It can be avoided if
we assume that
$\max_{1\leq i\leq m}\ep e^{t(\|\X^{i}_{1}\|+\|\Y^{i}_{1}\|)}\leq
\kappa$ for some $t>0$.

\begin{theorem}\label{th5} Under $H_{0}$, $d_{i}\equiv d$,
\textup{(C1$^{*}$)--(C3$^{*}$)}, we have for any $\epsilon>0$
%
\begin{eqnarray}\label{th5a}
&&\sup_{y\in R} \Bigl|\pr\Bigl(\max_{1\leq i\leq m}T^{2}_{ni}<
y \Bigr)-\exp\bigl(-mF_{n_{1},n_{2}}(y) \bigr) \Bigr|
\nonumber\\[-8pt]\\[-8pt]
&&\qquad \leq C \biggl(\sqrt{\frac{(\log m)^{5}}{n}}+m^{\rho
-(1-r)/(1+r)+\epsilon}+m^{\xi-1}\log
m \biggr),\nonumber
\end{eqnarray}
where $F_{n_1,n_2}(y)$ is defined in (\ref{a10}) and $C$ is a finite
constant depending only on $\xi,r,\rho,\delta,\kappa,\epsilon,c_{1},c_{2}$ and $d$.
\end{theorem}

If $m\geq c_{1}n^{b}$ for all $b>0$, then the error rate in Theorem
\ref{th5} is of order $\sqrt{(\log m)^{5}/n}$.
By Theorem~\ref{th5}, we can get the following result.

\begin{theorem}\label{th5-5} Under $H_{0}$, $d_{i}\equiv d$,
\textup{(C1$^{*}$)--(C3$^{*}$)}, we have for any $\epsilon>0$,
\[
\sup_{0\leq\alpha\leq1} \bigl|\pr\bigl(\Phi^{*}_{\alpha}=1 \bigr)-
\alpha\bigr|\leq C \biggl(\sqrt{\frac{(\log m)^{5}}{n}}+m^{\rho
-(1-r)/(1+r)+\epsilon}+m^{\xi-1}
\log m \biggr),
\]
where $C$ is given in (\ref{th5a}).
\end{theorem}

\subsection{\texorpdfstring{Power result for $\Phi_{\alpha}^{*}$}
{Power result for Phi alpha*}}\label{sec3.2}

Here we consider the power of the test $\Phi^{*}_{\alpha}$.

\begin{theorem}\label{th4} Suppose that
\[
\max_{1\leq i\leq m}\bigl\|\S^{-1/2}_{i}(\m_{1i}-
\m_{2i})\bigr\|\geq\sqrt{\frac{(2+\epsilon)\log m}{n_{1}}}
\]
for some $\epsilon>0$.
Then under \textup{(C1)} and \textup{(C2)},
\[
\pr\bigl(\Phi^{*}_{\alpha}=1\bigr)\rightarrow1 \qquad\mbox{as $n
\rightarrow\infty$}.
\]
\end{theorem}

Theorem~\ref{th4} shows that, in order to reject the null hypothesis
correctly, we only require $\max_{1\leq i\leq m}\|\S^{-1/2}_{i}(\m
_{1i}-\m_{2i})\|\geq\sqrt{\frac{(2+\epsilon)\log m}{n_{1}}}$. The
optimality\vspace*{1pt} of this lower bound when $d=1$ can be found in
\citet{CaiLiuXia}. We believe this lower bound remains optimal for
$d\geq2$ under some regularity conditions.

\section{Numerical results}\label{sec4}

\subsection{Simulation}\label{sec4.1}

In this section, we examine the numerical performance of the proposed
tests $\Phi^{*}_{\alpha}$ with $d=3$. We first compare
$\Phi^{*}_{\alpha }$ with $\Phi_{\alpha}$ to see the improvement of the
intermediate approximation and then compare $\Phi^{*}_{\alpha}$ to the
higher criticism (HC$^{*}$) test procedure
[\citet{DonJin04,HalJin10}], the test procedure proposed by
\citet{CheQin10} (C-Q) and the univariate $t$-test procedure based
on $\max_{1\leq i\leq dm}t^{2}_{i}$ (U-T), where $t_i$ is the two
sample $t$-statistic based on the $i$th coordinates of the
observations. The higher criticism test statistic is defined as
\citet{HalJin10}
\[
\mathrm{HC}^{*}=\max_{j:1/q\leq p_{(j)}\leq1/2} \biggl\{\frac{\sqrt
{q}(j/q-p_{(j)})}{\sqrt{p_{(j)}(1-p_{(j)})}}
\biggr\},
\]
where $q=3m$, $p_{j}=\pr(|N(0,1)|\geq|t_{i}|)$ and $p_{(j)}$ is the
$j$th $p$-value after sorting in ascending order. There are also other
versions\vadjust{\goodbreak} of HC$^{*}$ statistics [\citet{DonJin04}].
They perform similarly in our numerical studies. The critical values
$\alpha_{n}$ with significance level 0.05 are taken to be the
solutions to
$\pr(\mathrm{HC}^{*}\geq\alpha_{n})=0.05$ under that $p_{j}$,
$1\leq
j\leq3m$, are i.i.d. uniform $(0,1)$ distributed random variables.

Let
\begin{eqnarray*}
\bigl(\bigl(\X^{1}\bigr)^{\prime},\ldots,\bigl(
\X^{m}\bigr)^{\prime}\bigr)&=&\bigl(Z^{1}_{1},\ldots,Z^{3m}_{1}\bigr)\times\S^{1/2},
\\
\bigl(\bigl(\Y^{1}\bigr)^{\prime},\ldots,\bigl(
\Y^{m}\bigr)^{\prime}\bigr)&=&\bigl(Z^{1}_{2},\ldots,Z^{3m}_{2}\bigr)\times\S^{1/2}
\end{eqnarray*}
be $3m$-dimensional random vectors with covariance matrix
$\S$, where $\{Z_{i}^{j}\}$ are i.i.d. random variables. We consider
four distributions of $Z_{i}^{j}$,
$N(0,1)$, $t(5)$, exponential distribution with parameter 1 (Exp(1)),
and Gamma distribution
with shape and scale parameters $(2,2)$ (Gamma$(2,2)$). The covariance
matrix $\S$ is taken to be:

\begin{longlist}[(3)]
\item[(1)]
$\S_{1}=(0.9^{|j-i|})$;

\item[(2)] $\S_{2}=(\sigma_{ij})$, where $\sigma_{ij}=\max\{
1-|j-i|/(0.1*(3m)),0\}$;

\item[(3)] $\S_{3}=(\sigma_{ij})$, where $\sigma_{ij}=\max\{
1-|j-i|/(0.8*(3m)),0\}$.
\end{longlist}


$\S_{1}$ is an approximately bandable matrix. $\S_{2}$ is a $0.3m$
sparse matrix which has $0.3m$ nonzero entries in each row.
In $\S_{3}$, the number of nonzero entries in each row is $2.4 m$ and
the dependence between the variables becomes stronger than that in $\S_{2}$.

The sample sizes $(n_{1},n_{2})$ are taken to be
$(6,12)$, $(12,24)$, $(24,48)$ and $m$ takes {values} $50,100,200,400$.
We carry out 5000 simulations to obtain the empirical sizes {with
nominal significance level 0.05.} The results for $\S=\S_{1}$ are
summarized in Table~\ref{tabl1}.
The simulation results when $\S$ takes the other covariance matrices
are stated in the supplement material [\citet{LiuSha}] due to limit
of space. We can see that the empirical sizes of $\Phi^{*}_{\alpha}$
and Chen and Qin's test are close to $0.05$.
$\Phi^{*}_{\alpha}$ still performs well when the dependence becomes
stronger ($\S=\S_{2}$ and $\S_{3}$). However,
the empirical sizes of $\Phi_{\alpha}$ suffer very serious
distortions. This indicates the intermediate approximation in Section
\ref{sec3}
gains a lot of improvement on the accuracy of controlling type I
errors. The test procedure $\Phi^{*}_{\alpha}$ is robust to the tails
of distributions and the dependence. On the other hand, the empirical
sizes of HC$^{*}$ are much larger than $0.05$.
This shows that HC$^{*}$ statistic may be not robust to the dependence.
We have also done additional simulations and found that, when the
variables are independent but not normally distributed, HC$^{*}$
statistic may suffer serious distortions from the nominal significance level.

\begin{table}
\def\arraystretch{0.9}
\caption{Comparison of empirical sizes with nominal significance
level 0.05 ($\S=\S_{1}$)}\label{tabl1}
\begin{tabular*}{\tablewidth}{@{\extracolsep{\fill}}l c c c
cccc@{}}
\hline
& \multicolumn{4}{c}{$\bolds{N(0,1)}$} & \multicolumn{3}{c@{}}{$\bolds{t(5)}$}
\\[-4pt]
& \multicolumn{4}{c}{\hrulefill} & \multicolumn{3}{c@{}}{\hrulefill} \\
& $\bolds{m  \setminus (n_{1},n_{2})}$ & $\bolds{(6,12)}$ &
$\bolds{(12,24)}$& $\bolds{(24,48)}$& $\bolds{(6,12)}$&
$\bolds{(12,24)}$ & $\bolds{(24,48)}$\\
\hline
\hphantom{0}50 & $\Phi^{*}_{\alpha}$ & \textbf{0.0516}& \textbf{0.0466} & \textbf{0.0430}
&\textbf{0.0412} &\textbf{0.0374} & \textbf{0.0404} \\
&$\Phi_{\alpha}$ & 0.8965 & 0.4760 & 0.2285 &0.8641 &0.4312 & 0.2078
\\
&HC$^{*}$ & 0.5986 & 0.4348 & 0.3514 &0.6028 &0.4438 & 0.3534 \\
&C-Q & 0.0634 & 0.0644 & 0.0632 &0.0646 &0.0660 & 0.0644 \\
[4pt]
100 & $\Phi^{*}_{\alpha}$ & \textbf{0.0558} & \textbf{0.0483} &\textbf{0.0508}
&\textbf{0.0423} &\textbf{0.0360} & \textbf{0.0442} \\
& $\Phi_{\alpha}$ & 0.9694 & 0.5799 & 0.2711 &0.9542 &0.5315 & 0.2364
\\
&HC$^{*}$ & 0.7584 & 0.5228 & 0.4260 &0.7460 &0.5334 & 0.4100 \\
&C-Q & 0.0606 & 0.0620 & 0.0626 &0.0642 &0.0614 & 0.0592 \\
[4pt]
200 &$\Phi^{*}_{\alpha}$ & \textbf{0.0602} & \textbf{0.0584} &\textbf{0.0515}
&\textbf{0.0464} &\textbf{0.0393} & \textbf{0.0420} \\
& $\Phi_{\alpha}$ & 0.9958 & 0.7045 & 0.3238 &0.9916 &0.6380 & 0.2783
\\
&HC$^{*}$ & 0.9072 & 0.6492 & 0.4920 &0.8986 &0.6438 & 0.4672 \\
&C-Q & 0.0624 & 0.0584 & 0.0600 &0.0566 &0.0570 & 0.0574 \\
[4pt]
400 &$\Phi^{*}_{\alpha}$ & \textbf{0.0636} & \textbf{0.0609} &\textbf{0.0495}
&\textbf{0.0464} &\textbf{0.0402} & \textbf{0.0406} \\
& $\Phi_{\alpha}$ & 1.0000 & 0.8198 & 0.3781 &0.9996 &0.7571 & 0.3253
\\
&HC$^{*}$ & 0.9840 & 0.7876 & 0.5660 &0.9814 &0.7820 & 0.5642 \\
&C-Q & 0.0552 & 0.0592 & 0.0604 &0.0508 &0.0580 & 0.0588 \\
\hline
& \multicolumn{4}{c}{\textbf{Exp(1)}}&\multicolumn{3}{c@{}}{$\bolds{\operatorname{Gamma}(2,2)}$} \\
\hline
\hphantom{0}50 & $\Phi^{*}_{\alpha}$ & \textbf{0.0355}& \textbf{0.0392} & \textbf{0.0450}
&\textbf{0.0403} &\textbf{0.0468} & \textbf{0.0451} \\
&$\Phi_{\alpha}$ & 0.8441 & 0.4294 & 0.2226 &0.8675 &0.4473 & 0.2291
\\
&HC$^{*}$ & 0.5950 & 0.4492 & 0.3584 &0.5924 &0.4370 & 0.3604 \\
&C-Q & 0.0628 & 0.0622 & 0.0688 &0.0580 &0.0728 & 0.0666 \\
[4pt]
100& $\Phi^{*}_{\alpha}$ & \textbf{0.0404} & \textbf{0.0372} &\textbf{0.0519}
&\textbf{0.0436} &\textbf{0.0414} & \textbf{0.0524 } \\
& $\Phi_{\alpha}$ & 0.9409 & 0.5230 & 0.2625 &0.9557 &0.5521 & 0.2725
\\
&HC$^{*}$ & 0.7502 & 0.5296 & 0.4188 &0.7640 &0.5352 & 0.4212 \\
&C-Q & 0.0620 & 0.0626 & 0.0644 &0.0664 &0.0582 & 0.0598 \\
[4pt]
200 &$\Phi^{*}_{\alpha}$ & \textbf{0.0408} & \textbf{0.0364} &\textbf{0.0498}
&\textbf{0.0481} &\textbf{0.0435} &\textbf{0.0551} \\
&$\Phi_{\alpha}$ & 0.9882 & 0.6355 & 0.3105 &0.9923 &0.6671 & 0.3196
\\
&HC$^{*}$ & 0.8910 & 0.6358 & 0.4806 &0.9042 &0.6538 & 0.5014 \\
&C-Q & 0.0602 & 0.0608 & 0.0630 &0.0570 &0.0556 & 0.0610 \\
[4pt]
400 &$\Phi^{*}_{\alpha}$ & \textbf{0.0460} & \textbf{0.0355} &\textbf{0.0517}
&\textbf{0.0478} &\textbf{0.0449} & \textbf{0.0529} \\
& $\Phi_{\alpha}$ & 0.9987 & 0.7430 & 0.3671 &0.9997 &0.7810 & 0.3693
\\
&HC$^{*}$ & 0.9766 & 0.7788 & 0.5768 &0.9838 &0.7916 & 0.5762 \\
&C-Q & 0.0570 & 0.0590 & 0.0568 &0.0518 &0.0544 & 0.0572 \\
\hline
\end{tabular*}
\end{table}

To evaluate the power, we consider both approximately sparse model and
dense model. Let $\m_{1i}=0$ for $1\leq i\leq m$. Set $\m=(\mu
_{1},\ldots,\mu_{3m})=\ep((\Y^{1})^{\prime},\ldots,\break(\Y^{m})^{\prime})$ and
$\sigma^{2}=\Var(Z^{1}_{1})$. Consider

\textit{Model} 1 (\textit{approximately sparse case}). Let $\mu
_{i}=(-0.2)^{i-1}\times2\sqrt{\sigma^{2}\log m/n_{2}}$ for $1\leq
i\leq3m$.\vadjust{\goodbreak}

\textit{Model} 2 (\textit{dense case}). Let $\mu_{i}=0.2(-1)^{i-1}\times2\sqrt
{\sigma^{2}\log m/n_{2}}$ for $1\leq\break i\leq3m$.

Because of
the serious distortion of empirical sizes of $\Phi_{\alpha}$ and
HC$^{*}$, we do not consider the power of $\Phi_{\alpha}$ and HC$^{*}$.
We only report the power results for the normal distributions due to
the high similarity of the results with other distributions.
The reject region for $\max_{1\leq i\leq dm}t^{2}_{i}$ is
$[y_{n}(\alpha),\infty)$ with $d=1$ in
$F_{n_{1},n_{2}}(y)$ and $y_{n}(\alpha)$ satisfying
\[
\exp\bigl(-3mF_{n_{1},n_{2}}\bigl(y_{n}(\alpha)\bigr) \bigr)=1-
\alpha.
\]
This gives a much more accurate approximation than the extreme
distribution (results will not be reported here).

\begin{table}
\def\arraystretch{0.9}
\caption{Comparison of empirical powers ($\S=\S_{1}$)}\label{tabl2}
\begin{tabular*}{\tablewidth}{@{\extracolsep{\fill}}l ccccccc@{}}
\hline
& \multicolumn{4}{c}{\textbf{Model 1}}&\multicolumn{3}{c@{}}{\textbf{Model 2}}
\\[-4pt]
& \multicolumn{4}{c}{\hrulefill} & \multicolumn{3}{c@{}}{\hrulefill} \\
& $\bolds{m  \setminus (n_{1},n_{2})}$& $\bolds{(6,12)}$ &$\bolds{(12,24)}$
& $\bolds{(24,48)}$&$\bolds{(6,12)}$&$\bolds{(12,24)}$ & $\bolds{(24,48)}$\\
\hline
\hphantom{0}50 & $\Phi^{*}_{\alpha}$ & \textbf{0.7343}& \textbf{0.9327} &
\textbf{0.9758}
& \textbf{0.9453}& \textbf{0.9959} & \textbf{0.9994} \\
&C-Q & 0.0755 & 0.0739 & 0.0755 & 0.1369 & 0.1343 & 0.1404 \\
&U-T & 0.0766 & 0.0938 & 0.1064 & 0.0901& 0.0890 & 0.0862 \\
[4pt]
100 & $\Phi^{*}_{\alpha}$ &\textbf{0.7489} &\textbf{0.9538} &\textbf{0.9880} &
\textbf{0.9943}& \textbf{1.0000} & \textbf{1.0000} \\
& C-Q & 0.0704 & 0.0733 & 0.0720 & 0.2201 & 0.2250 & 0.2295 \\
&U-T & 0.0713 & 0.1001 & 0.0921 & 0.1019& 0.1137 & 0.0875 \\
[4pt]
200 &$\Phi^{*}_{\alpha}$ & \textbf{0.7451} &\textbf{0.9635} &\textbf{0.9937} &
\textbf{0.9998}& \textbf{1.0000} & \textbf{1.0000} \\
& C-Q & 0.0761 & 0.0665 & 0.0705 & 0.4289& 0.4365 & 0.4303 \\
&U-T & 0.0719 & 0.1058 & 0.0945 & 0.1278& 0.1507 & 0.1160 \\
[4pt]
400 &$\Phi^{*}_{\alpha}$ & \textbf{0.7520} & \textbf{0.9696} &\textbf{0.9957}
& \textbf{1.000}\hphantom{\textbf{0}}& \textbf{1.0000} & \textbf{1.0000} \\
& C-Q & 0.0633 & 0.0634 & 0.0636 & 0.7701 & 0.7997 & 0.8007 \\
&U-T & 0.0703 & 0.1089 & 0.0951 & 0.1414& 0.2062 & 0.1467 \\
\hline
\end{tabular*}
\end{table}

In Table~\ref{tabl2}, we only state the results when $\S=\S_{1}$. The other
simulation results are given in the supplement material [\citet{LiuSha}]. Note that in model 1,
$n\|\m\|^{2}/m^{1/2}\rightarrow0$. The power of \citet{CheQin10} is
low, as shown in Table~\ref{tabl2}. The power of $\max_{1\leq i\leq dm}t^{2}_{i}$
is also quite low. Our test statistics $\Phi^{*}_{\alpha}$ has the
highest powers which are close to one for $(n_{1},n_{2})=(12,24)$ and $(24,48)$.
In the dense case model 2, our test statistics still has the highest power.
We should remark that no method can uniformly outperform others over
all models and there may exist certain situations where Chen and Qin's
(\citeyear{CheQin10}) test statistic may outperform ours.

\subsection{Real data analysis}\label{sec4.2}

We apply the test procedure in Section~\ref{sec3} to test whether the tamoxifen
therapy is effective
on the promoter DNA methylation
status of 117 genes. 
The dataset consists of 123 patients, who showed the
extreme types of response to tamoxifen treatment; they either had an
objective response (CR${}+{}$PR, 45 patients) or a progressive disease
right from the start of treatment (PD, 78 patients). There are 117
genes and each gene
corresponds to a 2--6-dimensional vector that represents DNA
methylation status of CpG sites analyzed using a microarray-based DNA
methylation detection assay. \citet{Maretal05} used the
Benjamini--Hochberg (B-H) FDR procedure
with the target FDR of $25\%$ to identify
genes whose promoter DNA methylation
status was associated with the clinical benefit of tamoxifen therapy.
Before using B-H FDR\vadjust{\goodbreak} procedure, it is interesting to test whether
the tamoxifen therapy is effective on the promoter DNA methylation
status of those genes.

For each gene, we calculate the Hotelling's $T^{2}$-statistic
$T^{2}_{ni}$. The given significance level is $\alpha=0.05$. The value
of $\max_{1\leq i\leq m}G_{n_{1},n_{2},d_{i}}(T^{2}_{ni})$
is $1.0000$ which is larger than $1+m^{-1}\log(0.95)=0.9996$. Thus, we
can accept at the $0.05$ significance level that the tamoxifen therapy
has an effect on the promoter DNA methylation
status. We found three genes, PSAT1, STMN1 and SFN, whose values of
$G_{n_{1},n_{2},d_{i}}(T^{2}_{ni})$ are larger than $0.9996$. These
three genes were also identified by \citet{Maretal05} who used B-H
FDR correction and the $\chi^{2}$ distributions.\vspace*{-2pt}

\section{Proof of main results}\label{sec5}\vspace*{-2pt}

\subsection{\texorpdfstring{Proof of Theorem \protect\ref{th1}}{Proof of Theorem 2.1}}\label{sec5.1}

Without loss of generality, we assume that $\m_{1}=\m_{2}=0$. Since
$T_n^2$ converges to a chi-squared
distribution with $d$ degrees of freedom, we have for any $M>0$
\[
\lim_{n\rightarrow\infty}\sup_{0\leq x\leq M} \biggl|\frac{\pr
(T^{2}_{n}\geq x^{2})}{\pr(\chi^{2}(d)\geq x^{2})}-1 \biggr|=0.
\]
Thus, there exists a sequence $a_{n}\rightarrow\infty$ such that
%
\begin{equation}
\label{a0-0} \lim_{n\rightarrow\infty}\sup_{0\leq x\leq a_{n}} \biggl|\frac
{\pr
(T^{2}_{n}\geq x^{2})}{\pr(\chi^{2}(d)\geq x^{2})}-1 \biggr|=0.
\end{equation}
Let $\S=\S_{1}+\frac{n_{1}}{n_{2}}\S_{2}$ and
\[
\Z_{k}= \cases{\displaystyle  \S^{-1/2}\X_{k}, &\quad $1\leq k\leq
n_{1}$,
\vspace*{2pt}\cr
\displaystyle -\frac{n_{1}}{n_{2}}\S^{-1/2}\Y_{k-n_{1}}, &\quad
$n_{1}+1\leq k\leq n_{1}+n_{2}$.}
\]
By the identity
\[
\mathbf{x}^{\prime}\A^{-1}\mathbf{x}=
\max_{\|\theta\|=1}\frac
{(\mathbf{x}^{\prime}\theta)^{2}}{\theta^{\prime}\A\theta}
\]
for any $d\times d$ positive definite matrix $\A$, where $\theta$ is
a $d$-dimensional vector, we have
\begin{eqnarray*}
\bigl\{T^{2}_{n}\geq x^{2}\bigr\}&=& \Biggl\{
\exists\theta\mbox{, s.t. } \| \theta\|=1, \Biggl|\sum_{k=1}^{n}
\theta^{\prime}\Z_{k} \Biggr|
\\[-3pt]
&&\hspace*{7pt}\geq x\sqrt{\sum_{k=1}^{n}
\bigl(\theta^{\prime}\Z_{k}\bigr)^{2}-n_{1}
\bigl(\theta^{\prime}\bar{\Z}_{1}\bigr)^{2}-n_{2}
\bigl(\theta^{\prime}\bar{\Z}_{2}\bigr)^{2}} \Biggr\},
\end{eqnarray*}
where $n=n_{1}+n_{2}$, $\bar{\Z}_{1}=\frac{1}{n_{1}}\sum_{k=1}^{n_{1}}\Z
_{k}$ and $\bar{\Z}_{2}=\frac{1}{n_{2}}\sum_{k=n_{1}+1}^{n}\Z_{k}$.
Theorem~\ref{th1} follows if we can prove that
%
\begin{equation}
\label{a6} \frac{\pr( \exists\theta\mbox{, s.t. } \|\theta\|=1,
|\sum_{k\in H}\theta^{\prime}\Z_{k} |\geq
x\sqrt{\sum_{k\in H}(\theta^{\prime}\Z_{k})^{2}} )}{\pr(\chi^{2}(d)\geq
x^{2} )}\rightarrow1\vadjust{\goodbreak}
\end{equation}
uniformly for $x\in[a_{n}, o(n^{1/6}))$, $H=\{1,2,\ldots,n\}$, $\{
1,2,\ldots,n_{1}\}$
and $\{n_{1}+1,\ldots,n\}$. In fact, (\ref{a6}) implies that, for $i=1,2$,
\[
\frac{\pr(\exists\theta\mbox{, s.t. } \|\theta\|=1, |\theta^{\prime}\bar{\Z
}_{i}|\geq
2n^{-1}_{i}x\sqrt{\sum_{k=1}^{n}(\theta^{\prime}\Z_{k})^{2}} )}{\pr
(\chi^{2}(d)\geq4x^{2} )}\rightarrow1
\]
uniformly for $x\in[a_{n}, o(n^{1/6}))$. Observe that
\begin{eqnarray*}
\hspace*{-4pt}&&
\pr\bigl(T^{2}_{n}\geq x^{2}\bigr)
\\
\hspace*{-4pt}&&\qquad\leq \pr\Biggl( \exists\theta\mbox{, s.t. } \|\theta\|=1, \bigl|
\theta^{\prime}\bar{\Z}_{1}\bigr|\geq2n^{-1}_{1}x
\sqrt{\sum_{k=1}^{n_1}\bigl(
\theta^{\prime}\Z_{k}\bigr)^{2}} \Biggr)
\\
\hspace*{-4pt}&&\qquad\quad{} + \pr\Biggl( \exists\theta\mbox{, s.t. } \|\theta\|=1, \bigl|
\theta^{\prime}\bar{\Z}_{2}\bigr|\geq2n^{-1}_{2}x
\sqrt{\sum_{k=n_1+1}^{n}\bigl(
\theta^{\prime}\Z_{k}\bigr)^{2}} \Biggr)
\\
\hspace*{-4pt}&&\qquad\quad{} + \pr\biggl( \exists\theta\mbox{, s.t. } \|\theta\|=1,\frac{ |{\sum
_{k=1}^{n}}\theta^{\prime}\Z_{k} |
}{(\sum_{k=1}^n (\theta^{\prime} \Z_k)^2)^{1/2}}
\geq x \bigl(1-4x^{2}n^{-1}_{1}-4x^{2}n^{-1}_{2}
\bigr)^{1/2} \biggr)
\\
\hspace*{-4pt}&&\qquad = \bigl(2+o(1)\bigr) \pr\bigl(\chi^{2}(d)\geq4 x^{2}
\bigr)
\\
\hspace*{-4pt}&&\qquad\quad{} + \pr\biggl( \exists\theta\mbox{, s.t. } \|\theta\|=1, \frac{ |{\sum
_{k=1}^{n}}\theta^{\prime}\Z_{k} |
}{(\sum_{k=1}^n (\theta^{\prime} \Z_k)^2)^{1/2}}
\geq x \bigl(1-4x^{2}n^{-1}_{1}-4x^{2}n^{-1}_{2}
\bigr)^{1/2} \biggr)
\\
\hspace*{-4pt}&&\qquad = o(1) \pr\bigl(\chi^{2}(d)\geq x^{2} \bigr)
\\
\hspace*{-4pt}&&\qquad\quad{} + \pr\biggl( \exists\theta\mbox{, s.t. } \|\theta\|=1, \frac{ |{\sum
_{k=1}^{n}}\theta^{\prime}\Z_{k} |
}{(\sum_{k=1}^n (\theta^{\prime} \Z_k)^2)^{1/2}}
\geq x \bigl(1-4x^{2}n^{-1}_{1}-4x^{2}n^{-1}_{2}
\bigr)^{1/2} \biggr)
\end{eqnarray*}
uniformly in $x\in[a_{n}, o(n^{1/6}))$. Similarly, we can obtain a
lower bound for
$\pr(T^{2}_{n}\geq x^{2})$, which together with (\ref{a0-0}) and
(\ref{a6}) yields (\ref{th1a}).

We only prove (\ref{a6}) with $H=\{1,2,\ldots,n\}$. The proof for the
other two cases is similar.
Let $3/(3+\delta)<\beta<1$, $\hat{\Z}_{k}=\Z_{k}I\{\|\Z_{k}\|\leq
(\sqrt{n}/x)^{\beta}\}$ and set
\begin{eqnarray*}
S_{n}(\theta)&=&\sum_{k=1}^{n}
\theta^{\prime}\Z_{k},\qquad S^{\{\mathbf{N}\}}_{n}(\theta)=
\sum_{k=1, k
\notin
\mathbf{N}}^{n}\theta^{\prime}
\Z_{k},
\\
\hat{S}_{n}(\theta)&=&\sum_{k=1}^{n}
\theta^{\prime}\hat{\Z}_{k},\qquad \hat{S}^{\{\mathbf{N}\}}_{n}(
\theta)=\sum_{k=1, k\notin\mathbf{N}}^{n}
\theta^{\prime}\hat{\Z}_{k},
\\
\mathbf{V}_{n}(\theta)&=&\sum_{k=1}^{n}
\bigl(\theta^{\prime}\Z_{k}\bigr)^{2},\qquad
\mathbf{V}^{\{\mathbf{N}\}}_{n}(\theta)=\sum
_{k=1,k
\notin\mathbf{N}}^{n}\bigl(\theta^{\prime}
\Z_{k}\bigr)^{2},
\\
\hat{\mathbf{V}}_{n}(\theta)&=&\sum_{k=1}^{n}
\bigl(\theta^{\prime}\hat{\Z}_{k}\bigr)^{2},\qquad \hat{
\mathbf{V}}^{\{\mathbf{N}\}}_{n}(\theta)=\sum
_{k=1, k \notin
\mathbf{N}}^{n}\bigl(\theta^{\prime}\hat{
\Z}_{k}\bigr)^{2},
\end{eqnarray*}
where $\mathbf{N}$ is an index set.
By the fact that [see (5.7) in \citet{jing}]
%
\begin{equation}
\label{a1} \bigl\{s+t\geq x\sqrt{c+t^{2}}\bigr\}\subset\bigl\{s\geq
\bigl(x^{2}-1\bigr)^{1/2}\sqrt{c}\bigr\}
\end{equation}
for any $s,t\in R$, $c\geq0$ and $x\geq1$, we have
%
\begin{eqnarray}\label{a11}
&&\pr\bigl( \exists\theta\mbox{, s.t. } \|\theta\|=1,
\bigl|S_{n}(\theta)\bigr|\geq x\sqrt{ \mathbf{V}_{n}(\theta)}
\bigr)
\nonumber
\\
&&\qquad \leq \pr\bigl( \exists\theta\mbox{, s.t. } \|\theta\|=1, \bigl|
\hat{S}_{n}(\theta)\bigr|\geq x\sqrt{ \hat{\mathbf{V}}_{n}(
\theta)} \bigr)
\nonumber
\\
&&\qquad\quad{} + \sum_{j=1}^{n}\pr\bigl( \exists
\theta \mbox{, s.t. } \|\theta\|=1, \bigl|S^{\{j\}}_{n}(\theta)\bigr|\geq
\sqrt{x^{2}-1}\sqrt{ \mathbf{V}^{\{j\}}_{n}(
\theta)},A_{j} \bigr)
\\
&&\qquad= \pr\bigl(\exists\theta\mbox{, s.t. } \|\theta\|=1, \bigl|\hat
{S}_{n}(\theta)\bigr|\geq x\sqrt{ \hat{\mathbf{V}}_{n}(\theta)}
\bigr)
\nonumber
\\
&&\qquad\quad{} + \sum_{j=1}^{n}\pr\bigl( \exists
\theta\mbox{, s.t. } \|\theta\|=1, \bigl|S^{\{j\}}_{n}(\theta)\bigr|\geq
\sqrt{x^{2}-1}\sqrt{ \mathbf{V}^{\{j\}}_{n}(
\theta)} \bigr) P(A_{j}),\nonumber
\end{eqnarray}
where
\[
A_{j}=\bigl\{\|\Z_{j}\|\geq(\sqrt{n}/x)^{\beta}
\bigr\} \qquad\mbox{for $1\leq j\leq n$.}
\]
Repeating (\ref{a11}) and inequality (\ref{a1}) $m$ times, we get
\begin{eqnarray*}
&&
\pr\bigl(\exists\theta\mbox{, s.t. } \|\theta\|=1, \bigl|S_{n}(
\theta)\bigr|\geq x\sqrt{ \mathbf{V}_{n}(\theta)} \bigr)
\\
&&\qquad\leq\pr\bigl(\exists\theta\mbox{, s.t. } \|\theta\|=1, \bigl|\hat
{S}_{n}(\theta)\bigr|\geq x\sqrt{ \hat{\mathbf{V}}_{n}(\theta)}
\bigr)+\sum_{l=1}^{m}\hat{U}_{l}+U_{m+1},
\end{eqnarray*}
where
\begin{eqnarray*}
\hat{U}_{l}&=&\sum_{j_{1}=1}^{n}\cdots\sum_{j_{l}=1}^{n} \Biggl[\prod
_{k=1}^{l}\pr(A_{j_{k}}) \Biggr]
\\
&&\hspace*{49.5pt}{} \times\pr\bigl(\exists\theta\mbox{, s.t. } \|\theta\|=1, \bigl|
\hat{S}^{\{j_{1},\ldots,j_{l}\}}_{n}(\theta)\bigr| \geq\sqrt{x^{2}-l}\sqrt
{ \hat{\mathbf{V}}^{\{j_{1},\ldots,j_{l}\}}_{n}(\theta)} \bigr)
\end{eqnarray*}
and
\[
U_{m+1}=\sum_{j_{1}=1}^{n}\cdots\sum
_{j_{m+1}=1}^{n}\prod
_{k=1}^{m+1}\pr(A_{j_{k}}).
\]
Let $m=[x^{2}/2]$ for $x\geq4$. We have
%
\begin{eqnarray}
\label{ad1} U_{m+1}&=& \Biggl(\sum_{k=1}^{n}
\pr\bigl(\|\Z_{k}\|\geq(\sqrt{n}/x)^{\beta}\bigr)
\Biggr)^{m+1}
\nonumber\\[-8pt]\\[-8pt]
&\leq& e^{-m\log q_{n}}=o(1)\pr\bigl(\chi^{2}(d)\geq x \bigr),
\nonumber
\end{eqnarray}
where
\[
q_{n}= \bigl(n(x/\sqrt{n})^{\beta(3+\delta)}\ep\bigl(\|\X_{1}
\|^{3+\delta}+\|\Y_{1}\|^{3+\delta}\bigr) \bigr)^{-1}
\rightarrow\infty.
\]
The proof of (\ref{a6}) now relies on the Cram\'er-type moderate
theorem for
self-normalized truncated variables given below.

\begin{proposition} \label{prop1} Assume that $\operatorname{Card}(\mathbf{N})=O(x^{2})$. Then we have
%
\begin{eqnarray}
\label{prop1a} &&\pr\bigl(\exists\theta\mbox{, s.t. } \|\theta\|=1, \bigl|
\hat{S}^{\{
\mathbf{N}\}}_{n}(\theta)\bigr| \geq x\sqrt{ \hat{
\mathbf{V}}^{\{\mathbf{N}\}}_{n}(\theta)} \bigr)
\nonumber\\[-8pt]\\[-8pt]
&&\qquad =\bigl(1+o(1)\bigr)\pr\bigl(\chi^{2}(d)\geq x^{2} \bigr)
\nonumber
\end{eqnarray}
uniformly in $x\in[a_{n}, o(n^{1/6}))$.
\end{proposition}

The proof of Proposition~\ref{prop1} will be given in the next
subsection. Let us now finish the proof of
(\ref{a6}).

Using the same arguments as in the proof of inequality (\ref{ad1}) and
by Proposition~\ref{prop1}, we have
\begin{eqnarray*}
\sum_{l=1}^{m}\hat{U}_{l}&
\leq& C\sum_{l=1}^{m}\pr\bigl(
\chi^{2}(d)\geq x^{2}-l \bigr)\exp(-l\log q_{n})
\\
&=&o(1)\pr\bigl(\chi^{2}(d)\geq x^{2} \bigr)
\end{eqnarray*}
uniformly in $x\in[a_{n}, o(n^{1/6}))$. Hence,
\[
\pr\bigl(\exists\theta\mbox{, s.t. } \|\theta\|=1, \bigl|S_{n}(\theta)\bigr|
\geq x\sqrt{ \mathbf{V}_{n}(\theta)} \bigr)\leq\bigl(1+o(1)\bigr)\pr
\bigl(\chi^{2}(d)\geq x^{2} \bigr)
\]
uniformly in $x\in[a_{n}, o(n^{1/6}))$. To establish the lower bound,
we note that
\begin{eqnarray*}
&&\pr\bigl(\exists\theta\mbox{, s.t. } \|\theta\|=1, \bigl|S_{n}(
\theta)\bigr|\geq x\sqrt{ \mathbf{V}_{n}(\theta)} \bigr)
\\
&&\qquad\geq\pr\bigl(\exists\theta\mbox{, s.t. } \|\theta\|=1, \bigl|
\hat{S}_{n}(\theta)\bigr|\geq x\sqrt{ \hat{\mathbf{V}}_{n}(
\theta)} \bigr)
\\
&&\qquad\quad{} - \sum_{j=1}^{n}\pr\bigl(\exists
\theta\mbox{, s.t. } \|\theta\|=1, \bigl|\hat{S}^{\{j\}}_{n}(
\theta)\bigr|\geq\sqrt{x^{2}-1}\sqrt{ \hat{\mathbf{V}}^{\{j\}}_{n}(
\theta)} \bigr) P(A_{j}).
\end{eqnarray*}
It follows from Proposition~\ref{prop1} again that
\[
\pr\bigl(\exists\theta\mbox{, s.t. } \|\theta\|=1, \bigl|S_{n}(\theta)\bigr|
\geq x\sqrt{ \mathbf{V}_{n}(\theta)} \bigr)\geq\bigl(1+o(1)\bigr)\pr
\bigl(\chi^{2}(d)\geq x^{2} \bigr)
\]
uniformly in $x\in[a_{n}, o(n^{1/6}))$. This completes the proof of
(\ref{a6}) and hence Theorem~\ref{th1}.

\subsection{\texorpdfstring{Proof of Proposition \protect\ref{prop1}}{Proof of Proposition 5.1}}\label{sec5.2}

We start with the Cram\'er type moderate deviation theorem for
non-self-normalized sum.

\begin{lemma}\label{le2} Let $\operatorname{Card}(\mathbf{N})=O(x^{2})$.
We have
\[
\pr\bigl( \exists\theta\mbox{, s.t. } \|\theta\|=1, \bigl|\hat{S}^{\{
\mathbf{N}\}}_{n}(
\theta)\bigr|\geq x\sqrt{n_1} \bigr)=\bigl(1+o(1)\bigr)\pr\bigl(
\chi^{2}(d)\geq x^{2} \bigr)
\]
uniformly in $x\in[4, o(n^{1/6}))$.\vadjust{\goodbreak}
\end{lemma}

To prove Lemma~\ref{le2}, we need the following lemma by \citet{lin1}. The definition $|\cdot|_{d}$ below
is a slightly different from that in \citet{lin1}, but the proof
is exactly the same.

\begin{lemma}\label{le4} Let $\xi_{n,1},\ldots, \xi_{n,k_{n}}$ be
independent random vectors with mean zero
and values in $R^{d}$, and $S_{n}=\sum_{i=1}^{k_{n}} \xi_{n,i}$.
Assume that $\|\xi_{n,i}\|\leq
c_{n}B^{1/2}_{n}$, \mbox{$1\leq i\leq k_{n}$}, for some $c_{n}\rightarrow0$,
$B_{n}\rightarrow\infty$ and
\[
\bigl\|B^{-1}_{n}\Cov(\xi_{n,1}+\cdots+
\xi_{n,k_{n}})-I_{d} \bigr\| \leq C_{0}c^{2}_{n},
\]
where $I_{d}$ is a $d\times d$ identity matrix and $C_{0}$ is a
positive constant. Suppose that
$
\beta_{n}:=B^{-3/2}_{n}\sum_{i=1}^{k_{n}}\ep\|\xi_{n,i}\|^{3}\rightarrow
0.
$
Then for all $n\geq n_{0}$ ($n_{0}$ is given below)
\begin{eqnarray*}
&& \bigl|\pr\bigl(|S_{n}|_{d}\geq x\bigr)-\pr\bigl(|N|_{d}\geq
x/B^{1/2}_{n}\bigr)\bigr|
\\
&&\qquad\leq o(1)\pr\bigl(|N|_{d}\geq x/B^{1/2}_{n}
\bigr)\\
&&\qquad\quad{}+C_{d} \biggl(\exp\biggl(-\frac
{\delta^{2}_{n}\min(c^{-2}_{n}, \beta^{-2/3}_{n})}{8d} \biggr) +\exp
\biggl(\frac{C_{d}c^{2}_{n}}{\beta^{2}_{n}\log\beta_{n}} \biggr) \biggr),
\end{eqnarray*}
uniformly for $x\in[B^{1/2}_{n}, \delta_{n}\min(c^{-1}_{n}, \beta
^{-1/3}_{n})B^{1/2}_{n}]$, with any $\delta_{n}\rightarrow0$ and\break
$\delta_{n}\min(c^{-1}_{n}, \beta^{-1/3}_{n})\rightarrow\infty$,
where $N$ is a centered normal random vector
with covariance matrix $I_{d}$;
$|\cdot|_{d}$ denotes
$|\mathbf{z}|_{d}=\min\{\|\mathbf{x}_{i}\|\dvtx 1\leq i\leq d/q\}$,
$\mathbf{z}=(\mathbf{x}_{1},\ldots,\break\mathbf{x}_{d/q})$,
$\mathbf{x}_{i}\in R^{q}$ and $d/q$ is an integer; $o(1)$ is
bounded by $A_{n}:=A(\delta_{n}+\beta_{n})$, $A$ is
a positive constant depending only on $d$;
\[
n_{0}=\min\bigl\{n\dvtx  \forall k\geq n,c^{2}_{k}
\leq C_{01},\delta_{k}\leq C_{02},
\beta_{k}\leq C_{03} \bigr\},
\]
where $C_{01}$, $C_{02}$ and $C_{03}$ are some positive constants
depending only on $d$ and $C_{0}$.
\end{lemma}

\begin{pf*}{Proof of Lemma~\ref{le2}}
Let $\xi_{nk}=\hat{\Z}_{k}-\ep\hat {\Z}_{k}$, $B_{n}=n_{1}$ and
$c_{n}=\break 2n^{-1/2}_{1}(\sqrt{n}/ x)^{\beta }$ in Lemma~\ref{le4}. By the
inequalities $\beta>3/(3+\delta)$ and $x=o(n^{1/6})$,
\begin{eqnarray*}
\Biggl\|B_{n}^{-1}\Cov\Biggl(\sum_{k=1}^{n}
\xi_{nk}\Biggr)-I_{d} \Biggr\|&\leq& C\max_{1\leq k\leq n}\ep\|
\Z_{k}\|^{2}I\bigl\{\|Z_{k}\|\geq(\sqrt
{n}/x)^{\beta}\bigr\}
\\
&\leq&C(x/\sqrt{n})^{(1+\delta)\beta}\leq Cc_{n}^{2}.
\end{eqnarray*}
By letting $\delta_{n}\rightarrow0$ sufficiently slow, we have
\[
\exp\biggl(-\frac{\delta^{2}_{n}\min(c^{-2}_{n}, \beta
^{-2/3}_{n})}{8d} \biggr) +\exp\biggl(\frac{C_{d}c^{2}_{n}}{\beta
^{2}_{n}\log\beta_{n}}
\biggr)=o(1)\pr\bigl(\chi^{2}(d)\geq x^{2} \bigr)
\]
uniformly in $x\in[4, o(n^{1/6}))$. This proves Lemma~\ref{le2}.
\end{pf*}

\begin{pf*}{Proof of Proposition~\ref{prop1}} Observe that
\begin{eqnarray*}
&&\pr\bigl(\exists\theta\mbox{, s.t. } \|\theta\|=1, \bigl|\hat{S}^{\{
\mathbf{N}\}}_{n}(
\theta)\bigr|\geq x\sqrt{ \hat{\mathbf{V}}^{\{\mathbf{N}\}}_{n}(
\theta)} \bigr)
\\
&&\qquad\leq\pr\bigl(\exists\theta\mbox{, s.t. } \|\theta\|=1, \bigl|\hat
{S}^{\{\mathbf{N}\}}_{n}(\theta)\bigr|\geq x\sqrt{ n_{1}\bigl(1-
\varepsilon_{n}x^{-2}\bigr)} \bigr)
\\
&&\qquad\quad{} + \pr\bigl(\exists\theta\mbox{, s.t. } \|\theta\|=1, \bigl|\hat{S}^{\{
\mathbf{N}\}}_{n}(
\theta)\bigr|\geq x\sqrt{ \hat{\mathbf{V}}^{\{\mathbf{N}\}}_{n}(
\theta)},E_{n}(\theta) \bigr)
\end{eqnarray*}
and
\begin{eqnarray*}
&&\pr\bigl(\exists\theta\mbox{, s.t. } \|\theta\|=1, \bigl|\hat{S}^{\{
j\}}_{n}(
\theta)\bigr|\geq x\sqrt{ \hat{\mathbf{V}}^{\{j\}}_{n}(
\theta)} \bigr)
\\
&&\qquad\geq\pr\bigl(\exists\theta\mbox{, s.t. } \|\theta\|=1, \bigl|\hat
{S}^{\{j\}}_{n}(\theta)\bigr|\geq x\sqrt{ n_{1}\bigl(1+
\varepsilon_{n}x^{-2}\bigr)} \bigr)
\\
&&\qquad\quad{} - \pr\bigl(\exists\theta\mbox{, s.t. } \|\theta\|=1, \bigl|\hat{S}^{\{
j\}}_{n}(
\theta)\bigr|\geq x\sqrt{n_{1}\bigl(1+\varepsilon_{n}x^{-2}
\bigr)},F_{n}(\theta) \bigr),
\end{eqnarray*}
where $\varepsilon_{n}\rightarrow0$ which will be specified later and
\begin{eqnarray*}
E_{n}(\theta)&=&\bigl\{ \hat{\mathbf{V}}^{\{\mathbf{N}\}}_{n}(
\theta)\leq n_{1}\bigl(1-\varepsilon_{n}x^{-2}
\bigr)\bigr\},\\
F_{n}(\theta)&=&\bigl\{ \hat{\mathbf{V}}^{\{j\}}_{n}(
\theta)\geq n_{1}\bigl(1+\varepsilon_{n}x^{-2}
\bigr)\bigr\}.
\end{eqnarray*}
Also note that
\[
\pr\bigl( \exists\theta\mbox{, s.t. } \|\theta\|=1, \bigl|\hat{S}^{\{
\mathbf{N}\}}_{n}(
\theta)\bigr|\geq x\sqrt{n_{1}} \bigr)= \pr\bigl( \bigl|\hat
{S}^{\{\mathbf{N}\}}_{n}\bigr|_{d}\geq x\sqrt{n_{1}}
\bigr)
\]
with $q=d$.
By Lemma~\ref{le2}, we have
\[
\pr\bigl(\exists\theta\mbox{, s.t. } \|\theta\|=1, \bigl|\hat{S}^{\{
\mathbf{N}\}}_{n}(
\theta)\bigr|\geq x\sqrt{ n_{1}\bigl(1\pm\varepsilon_{n}x^{-2}
\bigr)} \bigr) =\bigl(1+o(1)\bigr)\pr\bigl(\chi^{2}(d)\geq
x^{2} \bigr)
\]
uniformly in $x\in[a_{n}, o(n^{1/6}))$. So it suffices to prove the
following lemma.
\end{pf*}

\begin{lemma}\label{a3} Let $\operatorname{Card}(\mathbf{N})=O(x^{2})$. We have
%
\begin{eqnarray}
\label{a4} &&\pr\bigl(\exists\theta\mbox{, s.t. }  \|\theta\|=1, \bigl|
\hat{S}^{\{
\mathbf{N}\}}_{n}(\theta)\bigr|\geq x\sqrt{ \hat{
\mathbf{V}}^{\{\mathbf{N}\}}_{n}(\theta)},E_{n}(\theta) \bigr)
\nonumber\\[-8pt]\\[-8pt]
&&\qquad =o(1)\pr\bigl(\chi^{2}(d)\geq x^{2} \bigr)
\nonumber
\end{eqnarray}
and
%
\begin{eqnarray}
\label{a5} &&\pr\bigl(\exists\theta\mbox{, s.t. } \|\theta\|=1, \bigl|
\hat{S}^{\{
j\}}_{n}(\theta)\bigr|\geq x\sqrt{n_{1}\bigl(1+
\varepsilon_{n}x^{-2}\bigr)},F_{n}(\theta) \bigr)
\nonumber\\[-8pt]\\[-8pt]
&&\qquad =o(1)\pr\bigl(\chi^{2}(d)\geq x^{2} \bigr)
\nonumber
\end{eqnarray}
uniformly in $x\in[a_{n}, o(n^{1/6}))$.
\end{lemma}

\begin{pf} We only prove (\ref{a4}) because the proof of
(\ref{a5}) is similar. Let $b=x/\sqrt{n_{1}}$. Then for
$0<\varepsilon_{n}<1/2$,
\begin{eqnarray*}
&&\bigl\{\hat{S}^{\{\mathbf{N}\}}_{n}(\theta)\geq x\sqrt{
\hat{\mathbf{V}}^{\{\mathbf{N}\}}_{n}(\theta)},E_{n}(\theta)
\bigr\}
\\
&&\quad \subset\bigl\{2b\hat{S}^{\{\mathbf{N}\}}_{n}(\theta)-b^{2}
\hat{\mathbf{V}}^{\{
\mathbf{N}\}}_{n}(\theta)\geq x^{2}-
\varepsilon^{2}_{n},E_{n}(\theta)\bigr\}
\\
&&\qquad{} \cup\bigl\{\hat{S}^{\{\mathbf{N}\}}_{n}(\theta)\geq x\sqrt
{ \hat{\mathbf{V}}^{\{\mathbf{N}\}}_{n}(\theta)},2xb\sqrt{ \hat
{\mathbf{V}}^{\{\mathbf{N}\}}_{n}(\theta)}<b^{2}\hat{\mathbf{V}}^{\{\mathbf{N}\}}_{n}(\theta)+x^{2}-\varepsilon^{2}_{n},E_{n}(
\theta)\bigr\}.
\end{eqnarray*}
We can choose $n_{d}$ points $\theta_{j}$, $1\leq j\leq n_{d}$,
with $\|\theta_{j}\|=1$ and $n_{d}\leq n^{2d}$, such that for any $\|
\theta\|=1$, $\|\theta-\theta_{j}\|\leq Cn^{-2}$ for some $1\leq
j\leq n_{d}$.
So we have
\begin{eqnarray*}
&&\pr\biggl(\bigcup_{\|\theta\|=1}\bigl\{2b
\hat{S}^{\{\mathbf{N}\}
}_{n}(\theta)-b^{2} \hat{
\mathbf{V}}^{\{\mathbf{N}\}}_{n}(\theta)\geq x^{2}-
\varepsilon^{2}_{n},E_{n}(\theta)\bigr\} \biggr)
\\
&&\qquad\leq\sum_{j=1}^{n_{d}}\pr\bigl(2b
\hat{S}^{\{\mathbf{N}\}
}_{n}(\theta_{j})-b^{2} \hat{
\mathbf{V}}^{\{\mathbf{N}\}
}_{n}(\theta_{j})\geq
x^{2}-\varepsilon^{2}_{n}-n^{-1}_{1},
\\
&&\hspace*{105.3pt} \mathbf{V}^{\{\mathbf{N}\}}_{n}(\theta_{j})\leq
n_{1}\bigl(1-\varepsilon_{n}x^{-2}
\bigr)+n^{-1}_{1} \bigr)
\\
&&\qquad\leq\sum_{j=1}^{n_{d}}\pr\bigl(2b
\hat{S}^{\{\mathbf{N}\}
}_{n}(\theta_{j})-b^{2}
\bigl(\hat{\mathbf{V}}^{\{\mathbf{N}\}
}_{n}(\theta_{j})-\ep
\hat{\mathbf{V}}^{\{\mathbf{N}\}}_{n}(\theta_{j})\bigr)\\
&&\qquad\quad\hspace*{25pt}{} + t
\bigl(\ep\hat{\mathbf{V}}^{\{\mathbf{N}\}}_{n}(\theta_{j})-
\hat{\mathbf{V}}^{\{\mathbf{N}\}}_{n}(\theta_{j})\bigr)
\\
&&\hspace*{38pt}\qquad \geq2x^{2}-\varepsilon^{2}_{n}-n^{-1}_{1}-O
\bigl(nb^{3}\bigr)+tn_{1}\varepsilon_{n}x^{-2}-O(ntb)
\bigr)
\\
&&\qquad=:\sum_{j=1}^{n_{d}}I_{j}.
\end{eqnarray*}
Let $t=(x/\sqrt{n})^{2-\gamma}$ with $0<\gamma<\beta(1+\delta)-1$
and $\max\{(x^{2}/n)^{\gamma/4},a^{-1/2}_{n}\}\leq\varepsilon
_{n}\rightarrow0$. We use Corollary 5 of \citet{Sak91} to bound
$I_{j}$. Let
\[
\xi_{k}=2b\theta^{\prime}_{j}\hat{
\Z}_{k}-2b\ep\theta^{\prime}_{j}\hat{\Z
}_{k}-\bigl(b^{2}-t\bigr) \bigl(\bigl(\theta^{\prime}_{j}
\hat{\Z}_{k}\bigr)^{2}-\ep\bigl(\theta^{\prime}_{j}
\hat{\Z}_{k}\bigr)^{2}\bigr),\qquad k\notin\mathbf{N}.
\]
Then $|\xi_{k}|=O(1)$, $B^{2}_{n}=\sum_{k\notin\mathbf{N}}\ep\xi
^{2}_{k}=4x^{2}+O(1)nb^{3}$, and for any bounded $h$,
\[
L(h)=\sum_{k\notin\mathbf{N}}\ep|\xi_{k}|^{3}
\max\bigl\{e^{h\xi
_{k}},1\bigr\}=O(1)nb^{3},
\]
where $O(1)$ are bounded by some absolute constants. Let
\[
y_{n}(x)=2x^{2}-\varepsilon^{2}_{n}-n^{-1}_{1}-O
\bigl(nb^{3}\bigr)+tn_{1}\varepsilon_{n}x^{-2}-O(ntb).
\]
By Corollary 5 of \citet{Sak91} and direct calculations, we obtain that
\begin{eqnarray*}
I_{j}&=&\bigl(1-\Phi\bigl(y_{n}(x)/B_{n}\bigr)
\bigr) \bigl(1+O\bigl(x^{3}/\sqrt{n}\bigr)\bigr)
\\
&=& O(1)x^{-1}\exp\bigl(-x^{2}/2-\bigl(n/x^{2}
\bigr)^{\gamma/2}\bigr)
\end{eqnarray*}
uniformly in $x\in[a_{n}, o(n^{1/6}))$. Hence, it follows that
%
\begin{eqnarray}
\label{a8} &&\pr\biggl(\bigcup_{\|\theta\|=1}\bigl\{2b
\hat{S}^{\{\mathbf{N}\}
}_{n}(\theta)-b^{2} \hat{
\mathbf{V}}^{\{\mathbf{N}\}}_{n}(\theta)\geq x^{2}-
\varepsilon^{2}_{n},E_{n}(\theta)\bigr\} \biggr)
\nonumber\\[-8pt]\\[-8pt]
&&\qquad=o(1)\pr\bigl(\chi^{2}(d)\geq x^{2} \bigr)
\nonumber
\end{eqnarray}
uniformly in $x\in[a_{n}, o(n^{1/6}))$.

Observe that
%
\begin{eqnarray}
\label{a7}\qquad &&\bigl\{\hat{S}^{\{\mathbf{N}\}}_{n}(\theta)\geq x\sqrt
{ \hat{\mathbf{V}}^{\{\mathbf{N}\}}_{n}(\theta)},2xb\sqrt
{ \hat{\mathbf{V}}^{\{
\mathbf{N}\}}_{n}(\theta)}<b^{2}
\hat{\mathbf{V}}^{\{\mathbf{N}\}
}_{n}(\theta)+x^{2}-
\varepsilon^{2}_{n},E_{n}(\theta)\bigr\}
\nonumber
\\
&&\qquad\subset\bigl\{\hat{S}^{\{\mathbf{N}\}}_{n}(\theta)\geq x\sqrt
{ \hat{\mathbf{V}}^{\{\mathbf{N}\}}_{n}(\theta)},b^{2}
\hat{\mathbf{V}}^{\{\mathbf{N}\}}_{n}(\theta)>x^{2}+
\varepsilon_{n} x,E_{n}(\theta)\bigr\}
\\
&&\qquad\quad{} \cup\bigl\{\hat{S}^{\{\mathbf{N}\}}_{n}(\theta)\geq x\sqrt{
\hat{\mathbf{V}}^{\{\mathbf{N}\}}_{n}(\theta)},b^{2}\hat{
\mathbf{V}}^{\{\mathbf{N}\}}_{n}(\theta)<x^{2}-
\varepsilon_{n} x,E_{n}(\theta)\bigr\}.\nonumber
\end{eqnarray}
By Lemma~\ref{le2},
\begin{eqnarray*}
&&\pr\biggl(\bigcup_{\|\theta\|=1}\bigl\{
\hat{S}^{\{\mathbf{N}\}
}_{n}(\theta)\geq x\sqrt{ \hat{
\mathbf{V}}^{\{\mathbf{N}\}
}_{n}(\theta)},b^{2}\hat{
\mathbf{V}}^{\{\mathbf{N}\}}_{n}(\theta)>x^{2}+
\varepsilon_{n} x,E_{n}(\theta)\bigr\} \biggr)
\\
&&\qquad\leq\pr\biggl(\bigcup_{\|\theta\|=1}\bigl\{
\hat{S}^{\{\mathbf{N}\}
}_{n}(\theta)\geq\sqrt{\bigl(x^{2}+
\varepsilon_{n} x\bigr)n_{1}}\bigr\} \biggr)
\\
&&\qquad=\bigl(1+o(1)\bigr)\pr\bigl(\chi^{2}(d)\geq x^{2}+
\varepsilon_{n} x \bigr)
\\
&&\qquad=o(1)\pr\bigl(\chi^{2}(d)\geq x^{2} \bigr)
\end{eqnarray*}
uniformly in $[a_{n},o(n^{1/6}))$ for any $a_{n}\rightarrow\infty$.
For the second term on the right-hand side of (\ref{a7}),
%
\begin{eqnarray}
\label{prop1-7}\quad &&\pr\biggl(\bigcup_{\|\theta\|=1}\bigl\{
\hat{S}^{\{\mathbf{N}\}
}_{n}(\theta)\geq x\sqrt{ \hat{
\mathbf{V}}^{\{\mathbf{N}\}
}_{n}(\theta)},b^{2}\hat{
\mathbf{V}}^{\{\mathbf{N}\}}_{n}(\theta)<x^{2}-
\varepsilon_{n} x,E_{n}(\theta)\bigr\} \biggr)
\nonumber
\\
&&\qquad\leq\sum_{k=1}^{[x]}\pr\biggl(\bigcup
_{\|\theta\|=1}\bigl\{\hat{S}^{\{
\mathbf{N}\}}_{n}(
\theta)\geq x\sqrt{\hat{\mathbf{V}}^{\{\mathbf{N}\}}_{n}(
\theta)},
\nonumber\\[-8pt]\\[-8pt]
&&\hspace*{93pt} \hat{\mathbf{V}}^{\{\mathbf{N}\}
}_{n}(\theta)\in
\bigl[n_{1}\bigl(1-\varepsilon_{n}(k+1)/x\bigr),
n_{1}(1-\varepsilon_{n}k/x)\bigr]\bigr\} \biggr)
\nonumber
\\
&&\qquad\quad{} +\pr\biggl(\bigcup_{\|\theta\|=1}\bigl\{\hat{
\mathbf{V}}^{\{\mathbf{N}\}}_{n}(\theta)\leq n_{1}(1-
\varepsilon_{n}/2)\bigr\} \biggr).\nonumber
\end{eqnarray}
For the last term above, we use the Bernstein inequality and obtain
\begin{eqnarray*}
&&\pr\biggl(\bigcup_{\|\theta\|=1}\bigl\{\hat{
\mathbf{V}}^{\{\mathbf{N}\}
}_{n}(\theta)\leq n_{1}(1-
\varepsilon_{n}/2)\bigr\} \biggr)
\\
&&\qquad\leq\sum_{j=1}^{n_{d}}\pr\bigl(\hat{
\mathbf{V}}^{\{\mathbf{N}\}
}_{n}(\theta_{j})\leq
n_{1}(1-\varepsilon_{n}/2)+n^{-1} \bigr)
\\
&&\qquad\leq\sum_{j=1}^{n_{d}}\pr\bigl(\ep\hat{
\mathbf{V}}^{\{\mathbf{N}\}}_{n}(\theta_{j})-\hat{
\mathbf{V}}^{\{\mathbf{N}\}}_{n}(\theta_{j}) \geq
n_{1}\bigl(\varepsilon_{n}/2+O(x/\sqrt{n})\bigr) \bigr)
\\
&&\qquad\leq\exp\biggl(-\frac{n_{1}(\varepsilon_{n}/2+O(x/\sqrt
{n}))^{2}}{2b^{-2\beta}+4b^{-2\beta}(\varepsilon_{n}/2+O(x/\sqrt
{n}))/3} \biggr)
\\
&&\qquad=o(1)\pr\bigl(\chi^{2}(d)\geq x^{2} \bigr)
\end{eqnarray*}
uniformly in $[a_{n},o(n^{1/6}))$. For the first term in (\ref{prop1-7}),
as in the proof of (\ref{a8}) using Corollary 5 of \citet{Sak91},
we can show that
\begin{eqnarray*}
&&\pr\biggl(\bigcup_{\|\theta\|=1}\bigl\{
\hat{S}^{\{\mathbf{N}\}
}_{n}(\theta)\geq x\sqrt{\hat{
\mathbf{V}}^{\{\mathbf{N}\}
}_{n}(\theta)},
\\[-2pt]
&&\hspace*{43.5pt} \hat{\mathbf{V}}^{\{\mathbf{N}\}
}_{n}(\theta)\in
\bigl[n_{1}\bigl(1-\varepsilon_{n}(k+1)/x\bigr),
n_{1}(1-\varepsilon_{n}k/x)\bigr]\bigr\} \biggr)
\\[-2pt]
&&\qquad\leq\pr\biggl(\bigcup_{\|\theta\|=1}\bigl\{
\hat{S}^{\{\mathbf{N}\}
}_{n}(\theta)\geq x\sqrt{n_{1}\bigl(1-
\varepsilon_{n}(k+1)/x\bigr)},
\\[-2pt]
&&\hspace*{116.5pt} \hat{\mathbf{V}}^{\{\mathbf{N}\}
}_{n}(\theta)\leq
n_{1}(1-\varepsilon_{n}k/x)\bigr\} \biggr)
\\[-2pt]
&&\qquad\leq\pr\biggl(\bigcup_{\|\theta\|=1}\bigl\{b
\hat{S}^{\{\mathbf{N}\}
}_{n}(\theta)+t\bigl(\ep\hat{\mathbf{V}}^{\{\mathbf{N}\}}_{n}(
\theta)-\hat{\mathbf{V}}^{\{\mathbf{N}\}}_{n}(\theta)\bigr)
\\[-2pt]
&&\hspace*{52pt}\qquad \geq x\sqrt{n_{1}\bigl(1-\varepsilon_{n}(k+1)/x\bigr)}
+n_{1}t\varepsilon_{n}k/x+O(ntb)\bigr\} \biggr)
\\[-2pt]
&&\qquad\leq Cn_{d}x^{-1}\exp\bigl(-x^{2}/2-c_{0}x^{-\gamma}n^{\gamma
/2}
\varepsilon_{n}\bigr)
\\[-2pt]
&&\qquad=o(1)\pr\bigl(\chi^{2}(d)\geq x^{2} \bigr)
\end{eqnarray*}
uniformly in $[a_{n},o(n^{1/6}))$. This completes the proof of Lemma
\ref{a3}.\vadjust{\goodbreak}
\end{pf}

\subsection{\texorpdfstring{Proof of Theorem \protect\ref{th3}}{Proof of Theorem 3.1}}\label{sec5.3}

Let $x_{n}=(2\log m+(d-2)\log\log m+x)^{1/2}$. Note that by Theorem
\ref{th1},
\[
\pr\Bigl(\max_{i\in\Lambda(r)}T^{2}_{ni}\geq
x^{2}_{n} \Bigr)\leq C\operatorname{Card}\bigl(\Lambda(r)
\bigr)m^{-1}=o(1).
\]
It suffices to prove that
\[
\pr\Bigl(\max_{i\notin\Lambda)(r)}T^{2}_{ni}\geq
x^{2}_{n} \Bigr)\rightarrow\exp\biggl(-\frac{1}{\Gamma(d/2)}
\exp(-x/2 ) \biggr).
\]
Since $\operatorname{Card}(\Lambda(r))=o(m)$, without loss of
generality, we
can assume that \mbox{$\Lambda(r)=\varnothing$},
that is, $\max_{1\leq i<j\leq m}\|\Gamma_{ij}\|\leq r$ for some
$r<1$. Otherwise, we only need to replace $\max_{1\leq i\leq m}(\cdot
)$ below by
$\max_{1\leq i\leq m,i\notin\Lambda(r)}(\cdot)$ and the proof
remains the same.
As in the proof of Theorem~\ref{th1}, we set
\[
\Z^{i}_{k}= \cases{ \S^{-1/2}_{i}
\X^{i}_{k}, &\quad $1\leq k\leq n_{1}$,
\vspace*{2pt}\cr
-
\displaystyle \frac{n_{1}}{n_{2}}\S^{-1/2}_{i}\Y^{i}_{k-n_1},
&\quad $n_{1}+1\leq k\leq n_{1}+n_{2}$,}
\]
and use the same truncation notations as in the proof of Theorem \ref
{th1}. With a careful check of the
proofs of Theorem~\ref{th1} and Proposition~\ref{prop1},
we can see that it suffices to show that, for $\operatorname{Card}(\mathbf{N})=O(x_{n}^{2})$,
%
\begin{eqnarray}
\label{th3-1}
\pr\Bigl(\max_{1\leq i\leq m} \bigl\|\hat{S}^{\{\mathbf{N}\}}_{ni}
\bigr\| \geq x_{n}\sqrt{ n_{1}\bigl(1\pm
\varepsilon_{n}x^{-2}_{n}\bigr)} \Bigr)
\rightarrow\exp\biggl(-\frac{1}{\Gamma(d/2)}\exp(-x/2 )
\biggr).\hspace*{-45pt}
\end{eqnarray}
Let $y_{n}=x_{n}\sqrt{ n_{1}(1\pm\varepsilon_{n}x^{-2}_{n})}$, where
$\varepsilon_{n}\rightarrow0$ to be specified
later. By the Bonferroni inequality, we have for any fixed integer $k$,
\begin{eqnarray*}
&&\sum_{l=1}^{2k}(-1)^{l-1}\sum
_{1\leq i_{1}<\cdots<i_{l}\leq m}\pr\bigl( \bigl\|\hat{S}^{\{\mathbf{N}\}}_{ni_{1}}
\bigr\|\geq y_{n},\ldots, \bigl\|\hat{S}^{\{\mathbf{N}\}}_{ni_{l}}\bigr\|\geq
y_{n} \bigr)
\\
&&\qquad\leq\pr\Bigl(\max_{1\leq i\leq m} \bigl\|\hat{S}^{\{\mathbf{N}\}
}_{ni}\bigr\|
\geq y_{n} \Bigr)
\\
&&\qquad\leq\sum_{l=1}^{2k-1}(-1)^{l-1}
\sum_{1\leq i_{1}<\cdots<i_{l}\leq
m}\pr\bigl( \bigl\|\hat{S}^{\{\mathbf{N}\}}_{ni_{1}}
\bigr\|\geq y_{n},\ldots, \bigl\|\hat{S}^{\{\mathbf{N}\}}_{ni_{l}}\bigr\|\geq
y_{n} \bigr).
\end{eqnarray*}
Theorem~\ref{th3} follows from the following lemma.

\begin{lemma}\label{le5} Let $\operatorname{Card}(\mathbf{N})=O(x^{2})$. We have for
any fixed $l$,
\begin{eqnarray*}
&&\sum_{1\leq i_{1}<\cdots<i_{l}\leq m}\pr\bigl( \bigl\|\hat{S}^{\{
\mathbf{N}\}}_{ni_{1}}
\bigr\|\geq y_{n},\ldots, \bigl\|\hat{S}^{\{\mathbf{N}\}}_{ni_{l}}\bigr\|\geq
y_{n} \bigr)
\\
&&\qquad =\bigl(1+o(1)\bigr)\frac{1}{l !} \biggl(\frac{1}{\Gamma(d/2)}\exp
(-x/2 )
\biggr)^{l}.
\end{eqnarray*}
\end{lemma}

In fact, by Lemma~\ref{le5}, we have
\begin{eqnarray*}
&&\limsup_{n\rightarrow\infty}\pr\Bigl(\max_{1\leq i\leq m} \bigl\|\hat
{S}^{\{\mathbf{N}\}}_{ni}\bigr\|\geq y_{n} \Bigr)
\\
&&\qquad\leq1-\sum_{l=0}^{2k-1}(-1)^{l}
\frac{1}{l !} \biggl(\frac{1}{\Gamma
(d/2)}\exp(-x/2 ) \biggr)^{l}
\\
&&\qquad\rightarrow1-\exp\biggl(-\frac{1}{\Gamma(d/2)}\exp(-x/2) \biggr)
\end{eqnarray*}
as $k\rightarrow\infty$. Similarly,
\[
\liminf_{n\rightarrow\infty}\pr\Bigl(\max_{1\leq i\leq m} \bigl\|\hat
{S}^{\{\mathbf{N}\}}_{ni}\bigr\|\geq y_{n} \Bigr)\geq1-\exp
\biggl(-\frac
{1}{\Gamma(d/2)}\exp(-x/2) \biggr).
\]
This proves Theorem~\ref{th3}.

\begin{pf*}{Proof of Lemma~\ref{le5}} Let $\X^{i}=(X^{i}_{1},\ldots,X^{i}_{d})^{\prime}$ and $\Y^{i}=(Y^{i}_{1},\ldots,Y^{i}_{d})^{\prime}$. Put
\[
r_{ij}=\max\Bigl\{\max_{k_{1},k_{2}}\bigl|\operatorname{Corr}
\bigl(X^{i}_{k_{1}}, X^{j}_{k_{2}}\bigr)\bigr|,
\max_{k_{1},k_{2}}\bigl| \operatorname{Corr}\bigl(Y^{i}_{k_{1}},
Y^{j}_{k_{2}}\bigr)\bigr|\Bigr\}
\]
and
\[
\mathcal{I}= \Bigl\{1\leq i_{1}<\cdots<i_{l}\leq m\dvtx
\max_{1\leq
k<j\leq l}r_{i_{k}i_{j}}\geq(\log m)^{-1-\gamma} \Bigr\}.
\]
When $l=1$, we let $\mathcal{I}=\varnothing$. For $2\leq j\leq l-1$, define
\begin{eqnarray*}
\mathcal{I}_{j}&=& \bigl\{1\leq i_{1}<\cdots<i_{l}\leq m\dvtx  \operatorname{Card}(\mathrm{S})=j\mbox{, where $
\mathrm{S}$ is the subset of }
\\
&&\hspace*{6pt}\{i_{1},\ldots,i_{l}\} \mbox{ with the largest
cardinality such that $\forall i_{k}\neq i_{t}\in
\mathrm{S}$,}
\\
&&\hspace*{188.5pt} r_{i_{k}i_{t}}< (\log m)^{-1-\gamma} \bigr\}.
\end{eqnarray*}
For $j=1$, define
\[
\mathcal{I}_{1} = \bigl\{1\leq i_{1}<\cdots<i_{l}\leq m\dvtx  r_{i_{k}i_{t}}\geq(\log m)^{-1-\gamma}
\mbox{ for every $1\leq k<t\leq l$} \bigr\}.
\]
It follows from the definition of $\mathcal{I}_{j}$ that $\mathcal
{I}=\bigcup_{j=1}^{l-1}\mathcal{I}_{j}$.
Then, by (C1), we have $\operatorname{Card}(\mathcal{I}_{j})=O(m^{j+2d\rho l})$. 
Define
\[
\mathcal{I}^{c}=\{1\leq i_{1}<\cdots<i_{l}\leq
m\}\setminus\mathcal{I}.
\]
We have $\operatorname{Card}(\mathcal{I}^{c})=C_{m}^{l}-O(m^{l-1+2d\rho
l})=(1+o(1))C_{m}^{l}$. For $(i_{1},\ldots,i_{l})\in\mathcal{I}^{c}$,
\[
\biggl\|\frac{1}{n_{1}}\Cov\bigl(\bigl(\hat{S}^{\{\mathbf{N}\}
}_{ni_{1}},\ldots,\hat{S}^{\{\mathbf{N}\}}_{ni_{l}}\bigr)\bigr)-I_{dl}\biggr\|\leq C(
\log m)^{-1-\gamma}+C(\log m/n)^{(1+\delta)\beta/2}.
\]
By Lemma~\ref{le4}, the proof of Lemma~\ref{le2} and some tedious
calculations,
\begin{eqnarray*}
&&\pr\bigl( \bigl\|\hat{S}^{\{\mathbf{N}\}}_{ni_{1}}\bigr\|\geq y_{n},\ldots, \bigl\|\hat{S}^{\{\mathbf{N}\}}_{ni_{l}}\bigr\|\geq y_{n} \bigr)
\\
&&\qquad=\bigl(1+o(1)\bigr)\pr\bigl( \|\W_{i_{1}}\|\geq y_{n}/
\sqrt{n_{1}},\ldots, \|\W_{i_{l}}\|\geq y_{n}/
\sqrt{n_{1}} \bigr),
\end{eqnarray*}
where $\W_{i_{1}},\ldots,\W_{i_{l}}$ are independent standard
$d$-dimensional random normal vectors.
By the tail probabilities of $\chi^{2}(d)$ distribution,
%
\begin{eqnarray}
\label{a9} &&\sum_{\mathcal{I}^{c}}\pr\bigl( \bigl\|
\hat{S}^{\{\mathbf{N}\}
}_{ni_{1}}\bigr\|\geq y_{n},\ldots, \bigl\|
\hat{S}^{\{\mathbf{N}\}}_{ni_{l}}\bigr\|\geq y_{n} \bigr)
\nonumber\\[-8pt]\\[-8pt]
&&\qquad=\bigl(1+o(1)\bigr)\frac{1}{l !} \biggl(\frac{1}{\Gamma(d/2)}\exp
(-y/2 )
\biggr)^{l}.
\nonumber
\end{eqnarray}
To prove the lemma, it suffices to show that for $1\leq j\leq l-1$,
%
\begin{equation}
\label{a10-0} \sum_{\mathcal{I}_{j}}\pr\bigl( \bigl\|
\hat{S}^{\{\mathbf{N}\}
}_{ni_{1}}\bigr\|\geq y_{n},\ldots, \bigl\|
\hat{S}^{\{\mathbf{N}\}}_{ni_{l}}\bigr\|\geq y_{n} \bigr)=o(1).
\end{equation}
To keep notation brief, we assume $\mathrm{S}=\{i_{l-j+1},\ldots,
i_{l}\}$ for $(i_{1},\ldots,i_{l})\in\mathcal{I}_{j}$. Divide
$\mathcal{I}_{j}$ into
$\mathcal{I}_{j1}$ and $\mathcal{I}_{j2}$, where
\begin{eqnarray*}
\mathcal{I}_{j1}&=&\biggl\{1\leq i_{1}<\cdots<i_{l}\leq m\mbox{: there
exists an $k\in\{i_{1},\ldots,i_{l-j}\}$}\\
&&\hspace*{8pt}\mbox{such that
for some $j_{1},j_{2}\in\mathrm{S}$ with
$\displaystyle j_{1}\neq j_{2}$, $\displaystyle r_{kj_{1}}\geq\frac{1}{(\log
m)^{1+\gamma}}$}\\
&&\hspace*{178pt}
\mbox{and $\displaystyle r_{kj_{2}}\geq\frac{1}{(\log m)^{1+\gamma}}$} \biggr\}
\end{eqnarray*}
and $\mathcal{I}_{j2}=\mathcal{I}_{j}\setminus\mathcal{I}_{j1}$.
Then $\operatorname{Card}(\mathcal{I}_{j1})=O(m^{j-1+4d\rho l})$ and again by Lemma
\ref{le4} and the proof of Lemma~\ref{le2},
\begin{eqnarray*}
&&\sum_{\mathcal{I}_{j1}}\pr\bigl( \bigl\|\hat{S}^{\{\mathbf{N}\}
}_{ni_{1}}
\bigr\|\geq y_{n},\ldots, \bigl\|\hat{S}^{\{\mathbf{N}\}}_{ni_{l}}\bigr\|\geq
y_{n} \bigr)
\\
&&\qquad \leq\sum_{\mathcal{I}_{j1}}\pr\bigl( \bigl\|\hat{S}^{\{\mathbf{N}\}
}_{ni_{l-j+1}}
\bigr\|\geq y_{n},\ldots, \bigl\|\hat{S}^{\{\mathbf{N}\}}_{ni_{l}}\bigr\|\geq
y_{n} \bigr)
\\
&&\qquad=\bigl(1+o(1)\bigr)\sum_{\mathcal{I}_{j1}}\pr\bigl( \|
\W_{i_{l-j+1}}\|\geq y_{n}/\sqrt{n_{1}},\ldots, \|
\W_{i_{l}}\|\geq y_{n}/\sqrt{n_{1}} \bigr)
\\
&&\qquad=O\bigl(m^{-1+4d\rho l}\bigr).
\end{eqnarray*}
For $(i_{1},\ldots,i_{l})\in\mathcal{I}_{j2}$ and $i_{l-j}$, there
is only one $j_{1}\in\mathrm{S}$
such that $r_{i_{l-j}j_{1}}\geq(\log m)^{-1-\gamma}$. For notation
briefness, we can assume $j_{1}=i_{l-j+1}$.
Thus, for any $0<\varepsilon<1$, by Theorem 1 in \citet{Zat87},
%
\begin{eqnarray}
\label{a14} &&\pr\bigl( \bigl\|\hat{S}^{\{\mathbf{N}\}}_{ni_{l-j}}\bigr\|\geq
y_{n},\ldots, \bigl\|\hat{S}^{\{\mathbf{N}\}}_{ni_{l}}\bigr\|\geq
y_{n} \bigr)
\nonumber
\\
&&\qquad\leq\pr\bigl( \|\tilde{\W}_{i_{l-j}}\|\geq(1-\varepsilon
)y_{n}/\sqrt{n_{1}},\ldots, \|\tilde{\W}_{i_{l}}\|
\geq(1-\varepsilon)y_{n}/\sqrt{n_{1}} \bigr)
\\
&&\qquad\quad{} +c_{1}\exp\bigl(-c_{2}(\log
m)^{1+(1-\beta)/2}\bigr),\nonumber
\end{eqnarray}
where $c_{1}$ and $c_{2}$ only depend on $d$ and $\varepsilon$,
$(\tilde{\W}_{i_{l-j}},\ldots,\tilde{\W}_{i_{l}})$ are
multivariate norm vector with covariance matrix $\Cov(\hat{S}^{\{
\mathbf{N}\}}_{ni_{l-j}},\ldots,\hat{S}^{\{\mathbf{N}\}
}_{ni_{l}})$. By the definition of $\mathcal{I}_{j2}$, we can prove that
\begin{eqnarray*}
&&
\biggl\|\frac{1}{n_{1}}\Cov\bigl(\hat{S}^{\{\mathbf{N}\}
}_{ni_{l-j}},\ldots,
\hat{S}^{\{\mathbf{N}\}}_{ni_{l}}\bigr)- \pmatrix{ \D & 0
\cr
0 & \mathbf{I} }
\biggr\|\\
&&\qquad\leq\frac{C}{(\log m)^{1+\gamma}}+C \biggl(\frac{\log
m}{n} \biggr)^{(1+\delta)\beta/2},
\end{eqnarray*}
where $\D=n_{1}^{-1}\sum_{k=1}^{n_{1}+n_{2}}\Cov((\Z^{i_{l-j}}_{k},\Z
_{k}^{i_{l-j+1}}))$ and $\mathbf{I}$ is
$(j-1)d$-dimensional identity matrix.
It follows that
\begin{eqnarray*}
&&\sum_{\mathcal{I}_{j2}}\pr\bigl( \bigl\|\hat{S}^{\{\mathbf{N}\}
}_{ni_{1}}
\bigr\|\geq y_{n},\ldots, \bigl\|\hat{S}^{\{\mathbf{N}\}}_{ni_{l}}\bigr\|\geq
y_{n} \bigr)
\\[-2pt]
&&\qquad \leq\sum_{\mathcal{I}_{j2}}\pr\bigl( \bigl\|\hat{S}^{\{\mathbf{N}\}
}_{ni_{l-j}}
\bigr\|\geq y_{n},\ldots, \bigl\|\hat{S}^{\{\mathbf{N}\}}_{ni_{l}}\bigr\|\geq
y_{n} \bigr)
\\[-2pt]
&&\qquad\leq\bigl(1+o(1)\bigr)\sum_{\mathcal{I}_{j2}}m^{-j+1}
\pr\bigl( \bigl\|(\tilde{\W}_{i_{l-j}},\tilde{\W}_{i_{l-j+1}})\bigr\|\geq(1-
\varepsilon)\sqrt{2}y_{n}/\sqrt{n_{1}} \bigr)\\[-2pt]
&&\qquad\quad{}+o(1).
\end{eqnarray*}
Since $\max_{1<i<j\leq p}\|\Gamma_{ij}\|\leq r$, we have $\|\D\|\leq
1+r$. This yields that
%
\begin{eqnarray}
\label{a13} &&\pr\bigl( \bigl\|(\tilde{\W}_{i_{l-j}},\tilde{\W}_{i_{l-j+1}})
\bigr\|\geq(1-\varepsilon)\sqrt{2}y_{n}/\sqrt{n_{1}} \bigr)
\nonumber\\[-8pt]\\[-8pt]
&&\qquad\leq C(\log m)^{d/2-1}m^{-2(1-\varepsilon)^{2}/(1+r)}.
\nonumber
\end{eqnarray}
Since $\rho$ is arbitrarily small, we can let $\varepsilon$ satisfy
$2(1-\varepsilon)^{2}/(1+r)> 1+\rho l$. This proves that
\[
\sum_{\mathcal{I}_{j2}}\pr\bigl( \bigl\|\hat{S}^{\{\mathbf{N}\}
}_{ni_{1}}
\bigr\|\geq y_{n},\ldots, \bigl\|\hat{S}^{\{\mathbf{N}\}}_{ni_{l}}\bigr\|\geq
y_{n} \bigr)=O\bigl(m^{j+\rho
l-j+1-2(1-\varepsilon)^{2}/(1+r)}\bigr)=o(1).
\]
Lemma~\ref{le5} is proved.
\end{pf*}

\subsection{\texorpdfstring{Proof of Theorem \protect\ref{th5}}{Proof of Theorem 3.3}}\label{sec5.4}

The proof of Theorem~\ref{th5} is given in the supplement material
[\citet{LiuSha}].

\subsection{\texorpdfstring{Proof of Theorem \protect\ref{th4}}{Proof of Theorem 3.5}}\label{sec5.5}

Let $i_{0}$ be the index such that
\[
\bigl\|\S^{-1/2}_{i_{0}}(\m_{1i_{0}}-\m_{2i_{0}})\bigr\|=
\max_{1\leq i\leq
m}\bigl\|\S^{-1/2}_{i}(\m_{1i}-
\m_{2i})\bigr\|\geq\sqrt{(2+\epsilon)\frac
{\log m}{n_{1}}}.
\]
Take\vspace*{1pt} $\|\theta\|=1$ such that $\theta^{\prime}\S^{-1/2}_{i_{0}}(\m
_{1i_{0}}-\m_{2i_{0}})=\|\S^{-1/2}_{i_{0}}(\m_{1i_{0}}-\m_{2i_{0}})\|$.
Note that $y_{n}(\alpha)=2\log m+(d-2)\log\log m+q_{\alpha}+o(1)$.
We have for any $0<\varepsilon<\sqrt{1+\epsilon/2}-1$,
\begin{eqnarray*}
\pr\bigl(\Phi^{*}_{\alpha}=1\bigr)&\geq& \pr
\bigl(T^{2}_{ni_{0}}\geq y_{n}(\alpha)\bigr)
\\[-2pt]
&\geq& \pr\Biggl( \sum_{k=1}^{n}
\theta^{\prime}\Z^{i_{0}}_{k}\geq(1+\varepsilon)
\sqrt{y_{n}(\alpha)n_{1}} \Biggr)+o(1)
\\[-2pt]
&\geq&\pr\Biggl( \sum_{k=1}^{n}
\theta^{\prime}\bigl(\Z^{i_{0}}_{k}-\ep
\Z^{i_{0}}_{k}\bigr)\geq(1+\varepsilon)\sqrt{y_{n}(
\alpha)n_{1}}-\sqrt{(2+\epsilon)n_{1}\log p} \Biggr)
\\[-2pt]
& &{}+o(1)
\\[-2pt]
&\rightarrow& 1.
\end{eqnarray*}

\section*{Acknowledgements}

We would like to thank the Associate Editor and referees for their insightful
comments that have led to significant improvement on the presentation
of the paper.

\begin{supplement}
\stitle{Supplement to ``A Cram\'{e}r moderate deviation theorem for
Hotelling's $T^{2}$-statistic with applications to global tests''}
\slink[doi]{10.1214/12-AOS1082SUPP} 
\sdatatype{.pdf}
\sfilename{aos1082\_supp.pdf}
\sdescription{The supplement material includes the moderate deviation result by
Sakhanenko (\citeyear{Sak91}), the proof of Theorem~\ref{th5} and the simulation results
in Section~\ref{sec4}.}
\end{supplement}


\printaddresses

\end{document}